\documentclass{svmult}
  
\usepackage{mathptmx,hyperref}   
\usepackage{helvet} 
\usepackage{courier}
\usepackage{amsmath,amsfonts}
\usepackage{graphicx}
\usepackage[bottom]{footmisc}          
    
\usepackage{mathptmx}         
\usepackage{helvet}            
\usepackage{courier}            
\usepackage{type1cm}             
\usepackage{makeidx}             
\usepackage{graphicx}          
\usepackage{multicol}          
\usepackage[bottom]{footmisc}  
     
\usepackage{epsfig}     
\usepackage{psfrag,rotating}     
\usepackage{amssymb,floatflt,enumerate} 
\usepackage{amsmath,amscd,psfrag,leqno} 
\usepackage[mathscr]{eucal}  
\usepackage{shadethm}           
\usepackage{graphicx,subfigure}   
\usepackage{overpic,contour}       
\contourlength{0.3mm}

\def\phm{{\hphantom{-}}}

\def\Vkt#1{{\mathbf #1}}

\makeindex

\begin{document}

\title*{On origami-like quasi-mechanisms with an antiprismatic skeleton}

\author{G. Nawratil}
\authorrunning{G. Nawratil}
\institute{
  Institute of Discrete Mathematics and Geometry \&  
	Center for Geometry and Computational Design, TU Wien,
  \email{nawratil@geometrie.tuwien.ac.at}}

%
%
\maketitle

\abstract{
We study snapping and shaky polyhedra which consist of antiprismatic skeletons covered by polyhedral belts composed of triangular faces only. 
In detail, we generalize Wunderlich's trisymmetric sandglass polyhedron in analogy to the 
generalizsation of the Jessen orthogonal icosahedron to Milka's extreme birosette structures, with 
the additional feature that the belt is developable into the plane as the Kresling pattern. 
Within the resulting $2$-dimensional family of origami-like sandglasses we study the 
1-parametric sets of quasi-mechanisms which are either {\it shaky} or have an {\it extremal snap}, i.e.\ one 
realization is on the boundary of self-intersection. 
Moreover, we evaluate the capability of these snapping/shaky quasi-mechanisms to flex
on base of the {\it snappability} index and the novel {\it shakeability} index, respectively.  
}


\section{Introduction}\label{sec:intro}

We consider polyhedral structures, where the bottom face $\alpha$ and the parallel top face $\beta$ are regular convex $n$-gons  $A_0,\ldots,A_{n-1}$ (with center $A$) and  $B_0,\ldots,B_{n-1}$ (with center $B$), respectively, 
for $n\geq 3$ and a side length of $1$. Moreover, these two faces are twisted against each other by a rotation about the orthogonal axis $AB$ 
(cf.\ Fig.\ \ref{fig1}a). Then this antiprismatic skeleton is covered by a polyhedral belt composed of triangular faces only. 
Moreover, all resulting polyhedra discussed in the paper do not possess continuous isometric deformations; i.e.\ they are rigid from the mathematical point of view.    
But the physical models of theses polyhedra can flex due to non-destructive elastic deformations of material (or backlash/tolerances in hinges), i.e.\ 
small changes in the intrinsic metric (given by the edge lengths) have significant effects on the spatial shape.  
Therefore these structures are also known as {\it model flexors} or {\it quasi-mechanisms}.

If the inner geometry (intrinsic metric together with the combinatorial structure) of the polyhedron is fixed, then the
embedding of the polyhedron into the Euclidean 3-space $E^3$ is in general not uniquely determined; i.e.\ 
the polyhedron has different incongruent realizations\footnote{In this paper the word realization always refers to an 
undeformed embedding into $E^3$.}. Based on these basic notations we can distinguish the following two
kinds of quasi-mechanisms:
\begin{enumerate}[(i)]
\item
Snapping quasi-mechanism: The shape variation results from the snap (caused by deformation) of a 
given realization into another one. An example for this is the Siamese dipyramid, which even snaps 
between three realizations \cite{goldberg,nawr2}.
\item
Shaky quasi-mechanism: Now the deformed states originate from a given {\it shaky} (also known as 
{\it singular} or {\it infinitesimal flexible}) realization. 
The best known example for this kind of model flexion is the Jessen orthogonal icosahedron \cite{jessen,goldberg,kalinin} (cf.\ Fig.\ \ref{fig2}a).
Note that shakiness can be seen as the limit of snapping when the related realizations converge to coincidence 
\cite{stachel_between,nawr2}.
\end{enumerate}
There are also examples, like the four-horn \cite{schwabe,nawr2},  which use both functional principles (i \& ii)
at the same \medskip time. 

\noindent
{\bf Review.} We proceed with a review on model flexors with an antiprismatic skeleton:

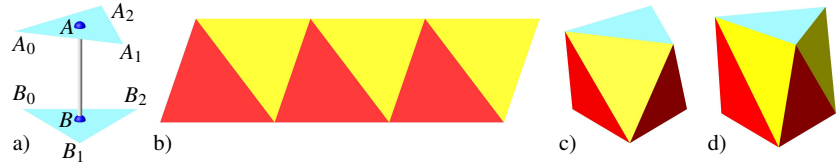
\begin{figure}[t]
\begin{center}
\begin{overpic}
    [width=110mm]{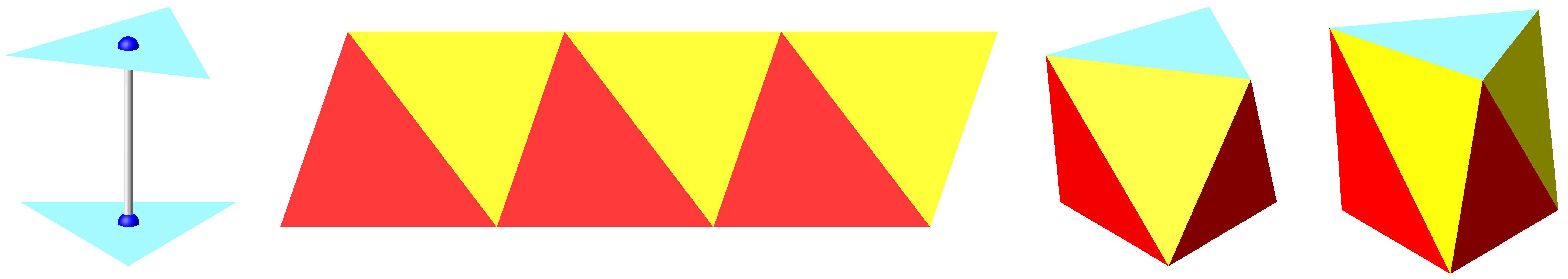}
\begin{small}
\put(0,0){a)}
\put(0,6){$B_0$}
\put(6,-1){$B_1$}
\put(13.5,5.8){$B_2$}
\put(5.6,2.8){$B$}
\put(0,12){$A_0$}
\put(13,11){$A_1$}
\put(12,16){$A_2$}
\put(5.8,14){$A$}
\put(17,0){b)}
\put(66,0){c)}
\put(84,0){d)}
\end{small}     
  \end{overpic}
\end{center}	
\caption{The antriprismatic skeleton for $n=3$ (a), its belt given by the Kresling pattern (b), the resulting octahedral structure (c) and its 
configuration after the snap (d).
}
  \label{fig1}
\end{figure}

\begin{enumerate}[(1)]
\item
Let us start with the Kresling pattern \cite{kresling}, a flat strip of congruent triangles (cf.\ Fig.\ \ref{fig1}b) which can be folded up and closed to a belt for the skeleton. 
The resulting antiprismatic structure has a bi-stable behavior, which was already known to  Wunderlich (cf.\ \cite{wunderlich_antiprism} and for $n=3$ 
the more detailed study \cite{wunderlich_achtflach}). 
During the snap, the relative motion of $\alpha$ and $\beta$ is composed of a rotation about $AB$ and a change in height (see Fig.\ \ref{fig1}c,b). 
A detailed literature review on these structures (as well as a study of related ones) is given by the author in \cite{nawr3}.
\item
Next, we consider so-called extreme birosettes \cite{birosettes}, which can be seen as 
generalizations of the Jessen orthogonal icosahedron (case $n=3$). 
 In this case the belt consists of $2n$ equilateral triangles with side length of $1$ and 
$2n$ {\it petals}, which are skew rhombi of side length 1 broken along one of its diagonals\footnote{Note that the birosette 
degenerates into an antiprism if $p$  converges towards zero.} of length $p$. 
The resulting structure is symmetric with respect to (a) rotations of $\tfrac{2\pi}{n}$ about $AB$ and (b) a reflection 
at the birosette center (midpoint of $AB$) which has to be complemented by a rotation of $\tfrac{\pi}{n}$ about $AB$ 
for even values of $n$ (see Figs.\ \ref{fig2}a and \ref{fig4}). The relative instantaneous motion of $\alpha$ and $\beta$ is translatory in direction of $AB$. 
Moreover, the maximal value $p_{+}$ of $p$ such that the birosette can be assembled under consideration of the mentioned symmetry, 
yields the {\it extreme} birosette. Note that birosettes for $n=3$ are also known as Douady shaddocks \cite{douady,kalinin}.
\item
Finally, we recall Wunderlich's trisymmetric sandglass polyhedron \cite{sanduhr}, which is  a snapping icosahedron possessing the same 
symmetry and combinatorial structure as the birosette for $n=3$. In this case the belt consists of six 
congruent isosceles triangles, whose bases of length $1$ are hinged to the equilateral triangles $\alpha$ and $\beta$, 
respectively, and the gaps between them are filled by further $12$ congruent isosceles triangles (see Fig.\ \ref{fig2}d,e). 
During the snap, the relative motion of $\alpha$ and $\beta$ is a translation along $AB$. 
\end{enumerate}

\begin{figure}[t]
\hfill
\begin{overpic}
    [width=20mm]{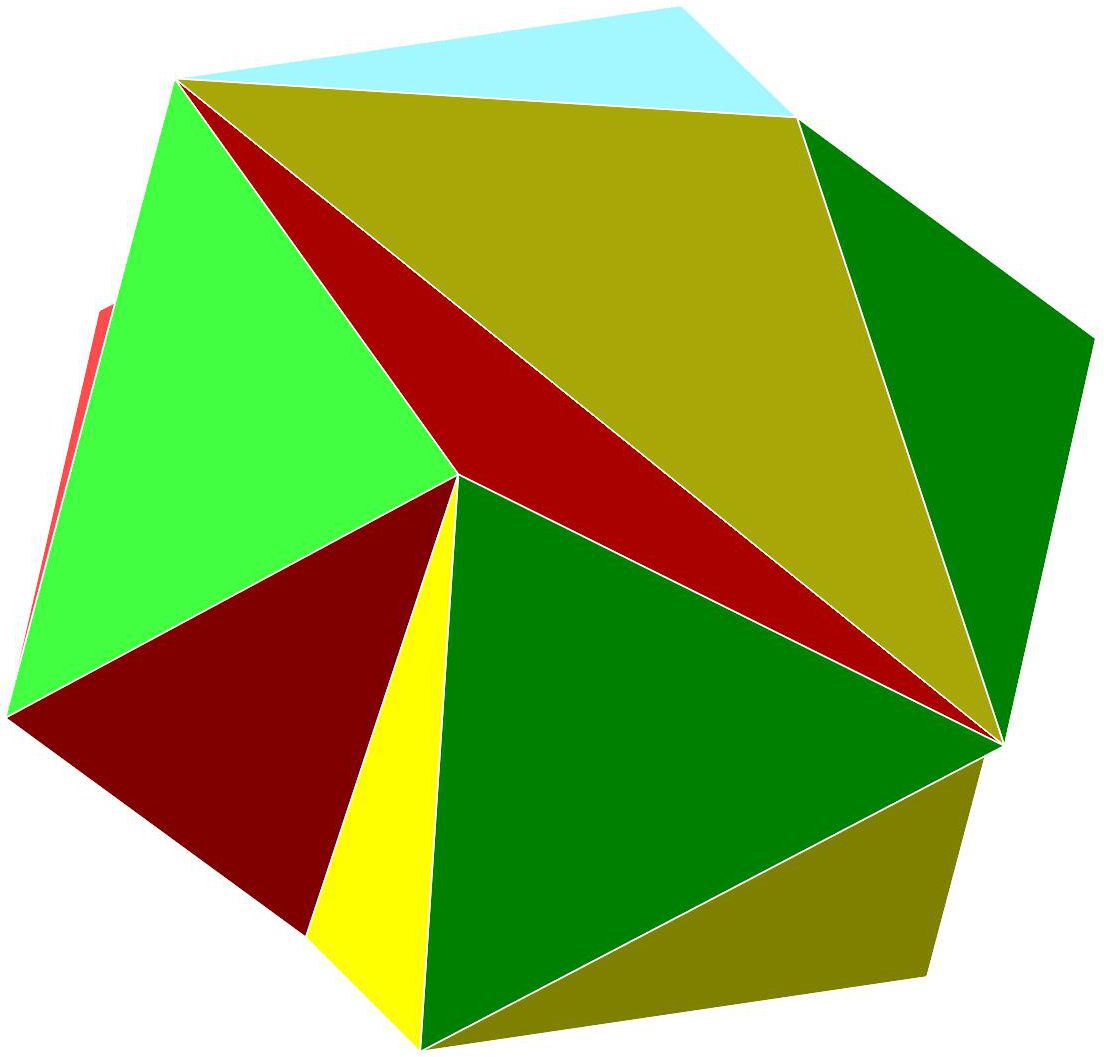}
\begin{small}
\put(0,0){a)}
\end{small}     
  \end{overpic} 
\hfill
\begin{overpic}
    [width=23mm]{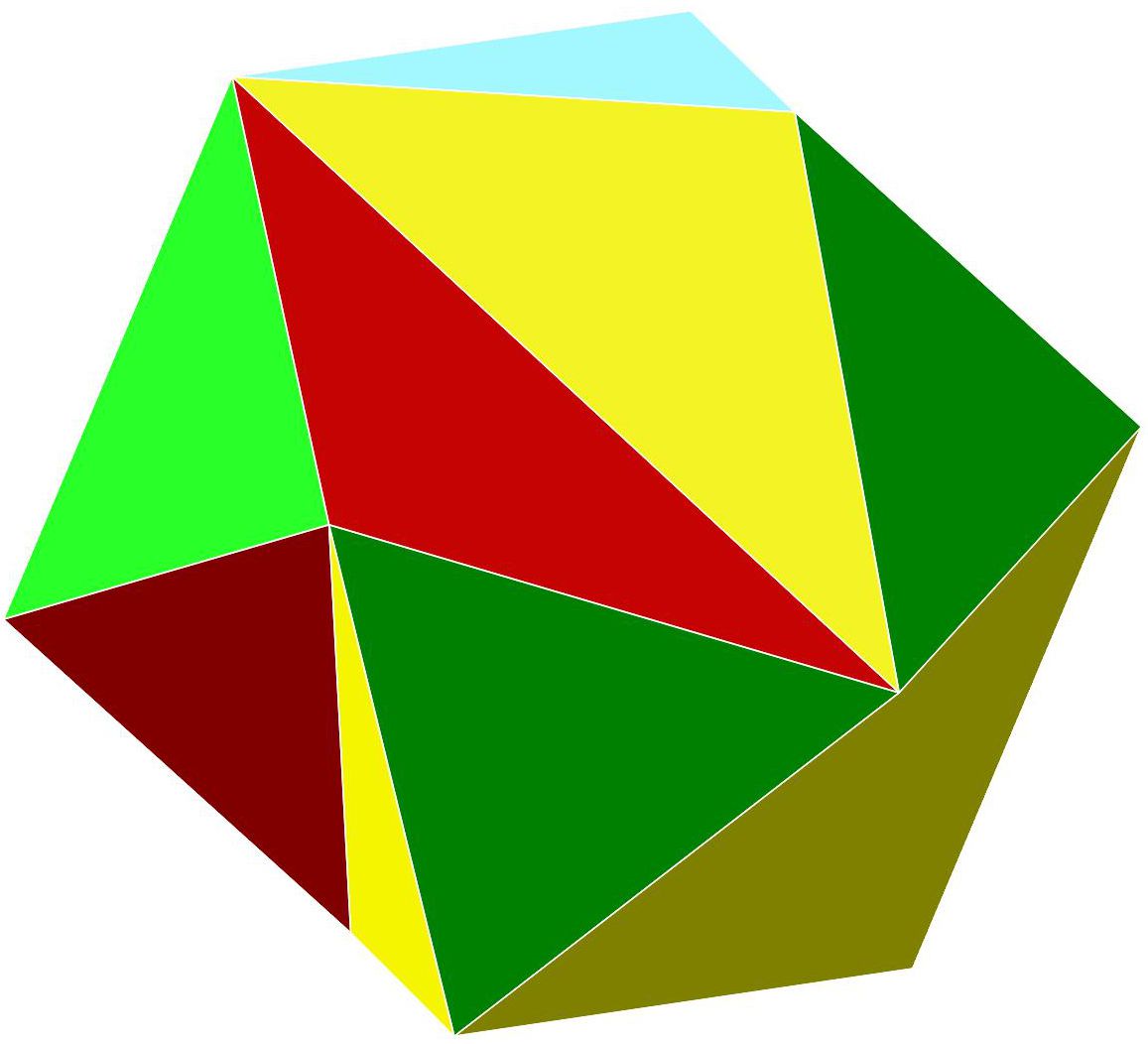}
\begin{small}
\put(0,0){b)}
\end{small}         
  \end{overpic} 
	\hfill
\begin{overpic}
    [width=17.5mm]{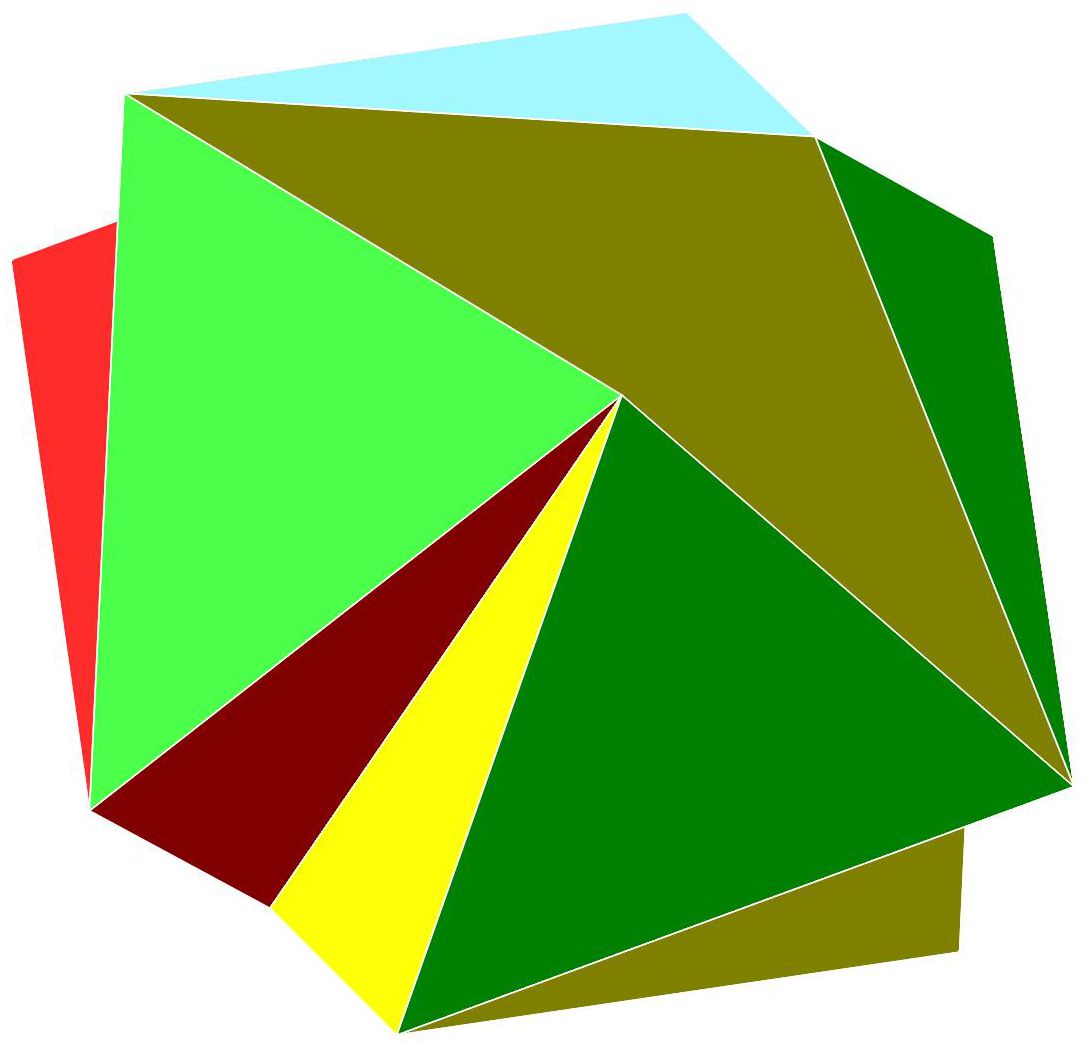}
\begin{small}
\put(0,0){c)}
\end{small}         
  \end{overpic} 	
\hfill
 \begin{overpic}
    [width=14mm]{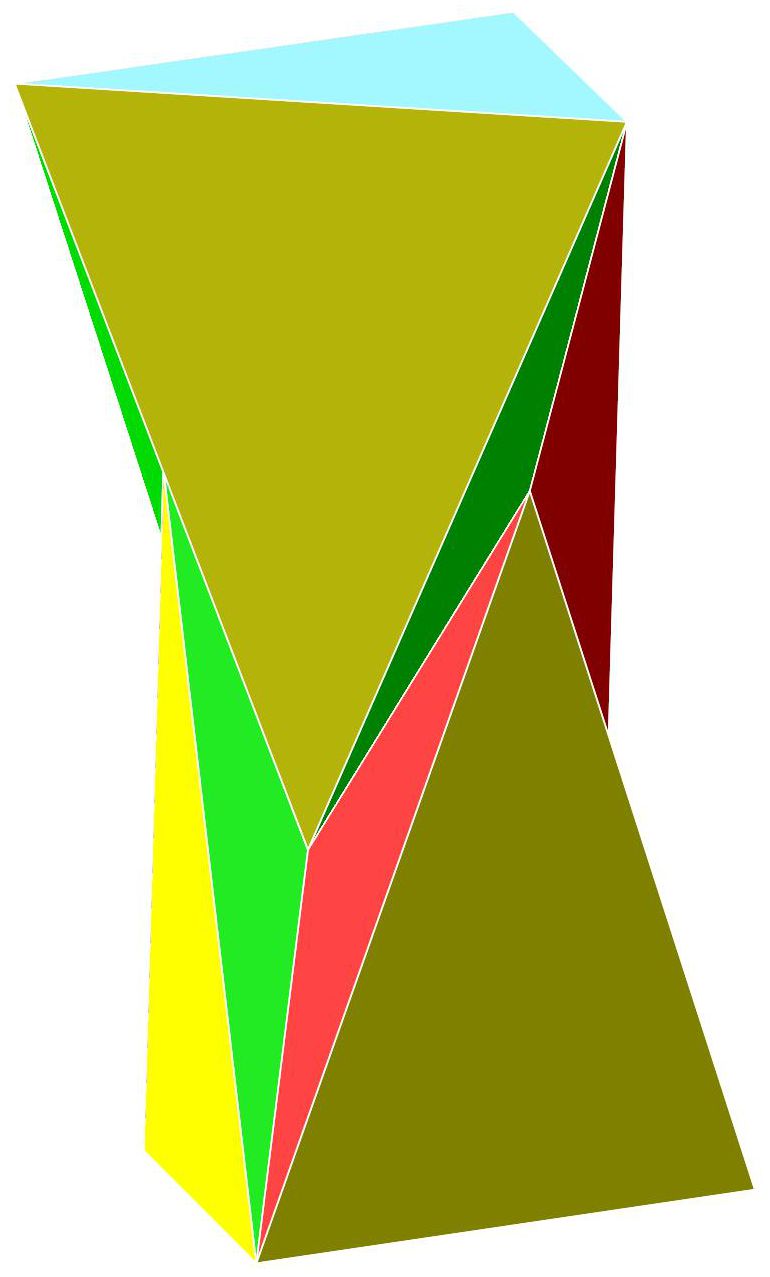}
\begin{small}
\put(0,0){d)}
\end{small}         
  \end{overpic} 
	\hfill
\begin{overpic}
    [width=14mm]{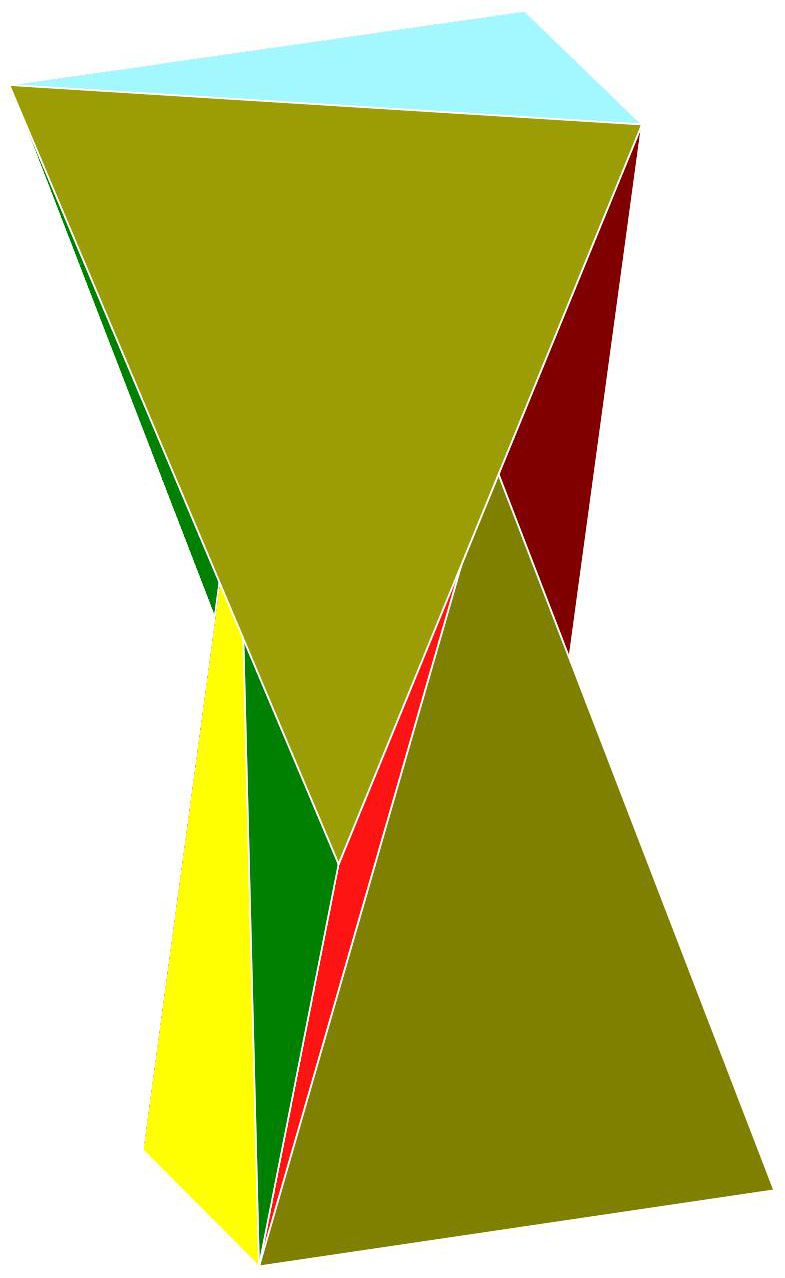}
\begin{small}
\put(0,0){e)}
\end{small}         
  \end{overpic} 	
\hfill
 \begin{overpic}
    [width=14mm]{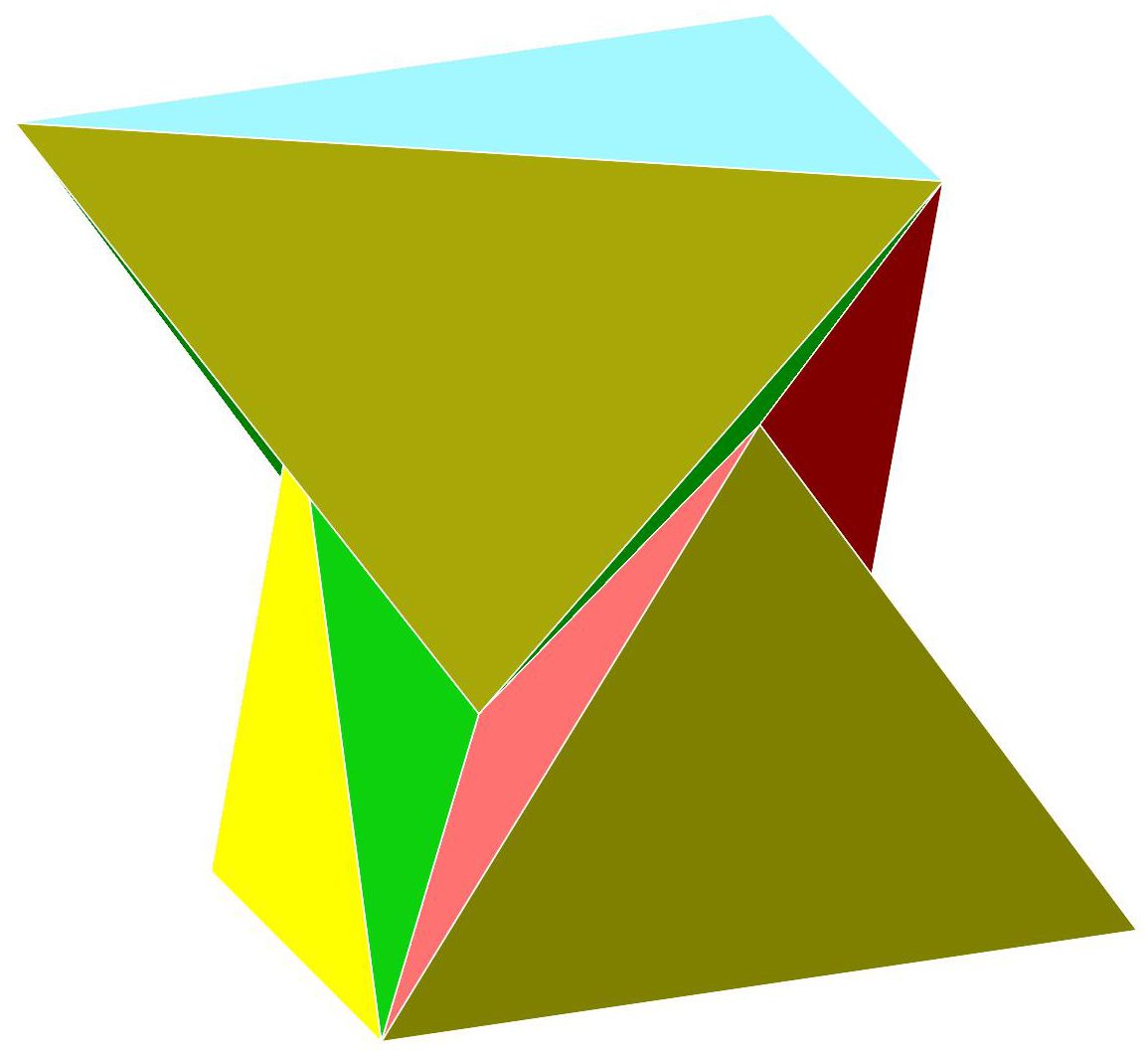}
\begin{small}
\put(0,0){f)}
\end{small}         
  \end{overpic} 
\caption{(a) Jessen orthogonal icosahedron \cite{jessen}, 
(b,c) two snapping realizations of a non-extremal birosette for $n=3$ with $p=1.6 < p_+=4/\sqrt{6}$,
(d,e) the two realizations of the snapping  sandglass icosahedron of Wunderlich \cite[Fig.\ 6]{sanduhr}, 
(f) shaky sandglass icosahedron \cite[Fig.\ 4]{sanduhr}. 
}
  \label{fig2}
\end{figure}

\begin{remark}
By considering the limit mentioned in (ii) one can also produce shaky quasi-mechanisms from the snapping structures given in (1) studied in \cite{wunderlich_antiprism,wunderlich_achtflach} 
and (3) illustrated in Fig.\ \ref{fig2}f. 
Also snapping birosettes can be generated by choosing $p<p_+$ as displayed in Fig.\ \ref{fig2}b,c. \hfill $\diamond$
\end{remark}

\noindent
{\bf Outline.}
Based on this review we generalize Wunderlich's sandglass polyhedron to arbitrarry $n$ in analogy to the birosette construction, with 
the additional feature that the belt is developable into the plane as the Kresling pattern (cf.\ Section \ref{sec:prelimi}), which 
allows an efficient production. 
Within the resulting $2$-dimensional family of origami-like sandglasses (for arbitrary $n\geq 3$) we study the 
1-parametric sets of quasi-mechanisms which are either {\it shaky} (cf.\ Section \ref{sec:shake}) or have an {\it extremal snap} (cf.\ Section \ref{sec:snap}), i.e.\ one 
realization is on the boundary of self-intersection. 
Moreover, in these two sections we also evaluate the capability of these snapping/shaky quasi-mechanisms to flex
on base of the {\it snappability} index (cf.\ \cite{nawr1,nawr2}) and the novel {\it shakeability} index, respectively.   
The paper is concluded in Section \ref{sec:end}.

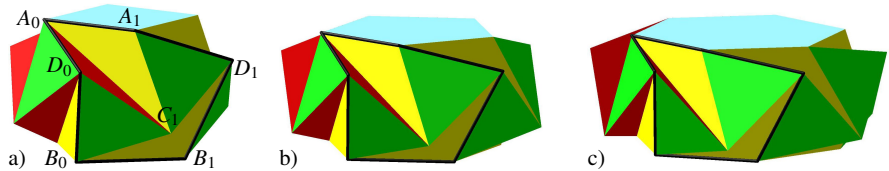
\begin{figure}[b]
\begin{overpic}
    [width=30mm]{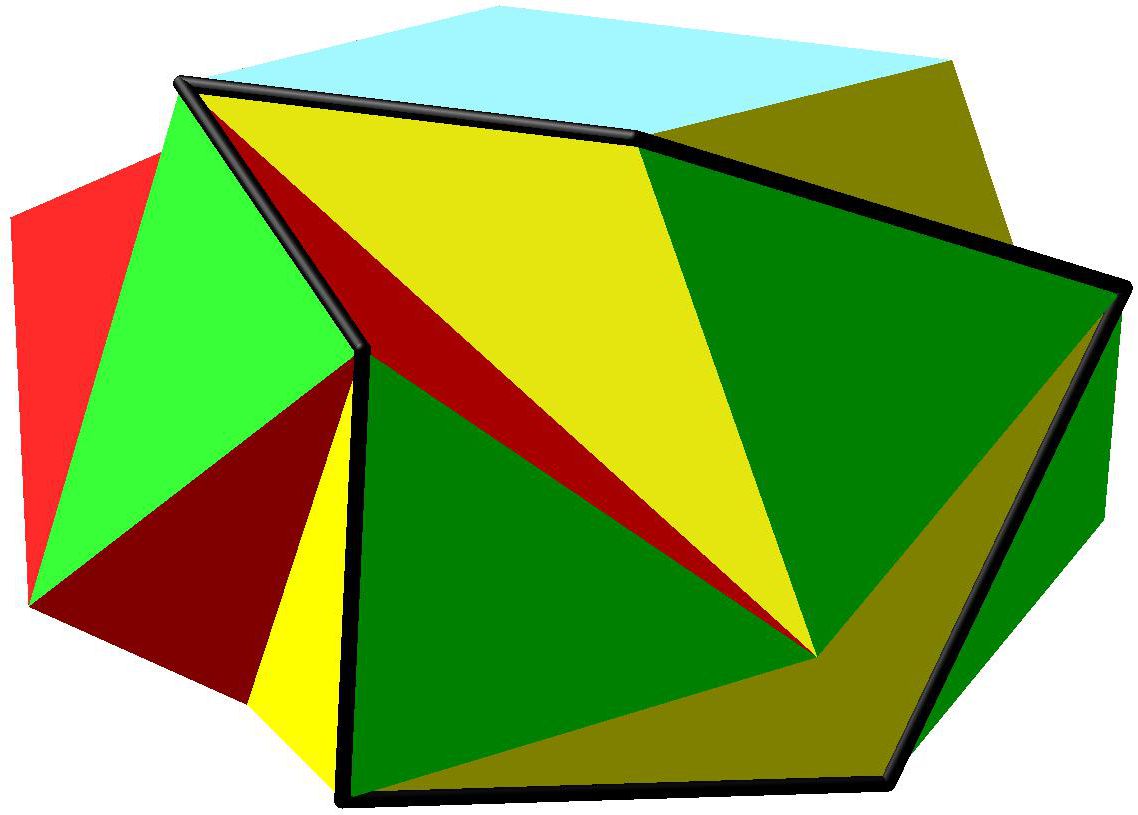}
\begin{small}
\put(0,0){a)}
\put(16,0){$B_0$}
\put(82,0){$B_1$}
\put(100,40){$D_1$}
\put(49,63){$A_1$}
\put(4,62){$A_0$}
\put(17,41){$D_0$}
\put(66,19){$C_1$}
\end{small}     
  \end{overpic} 
\hfill
 \begin{overpic}
    [width=35mm]{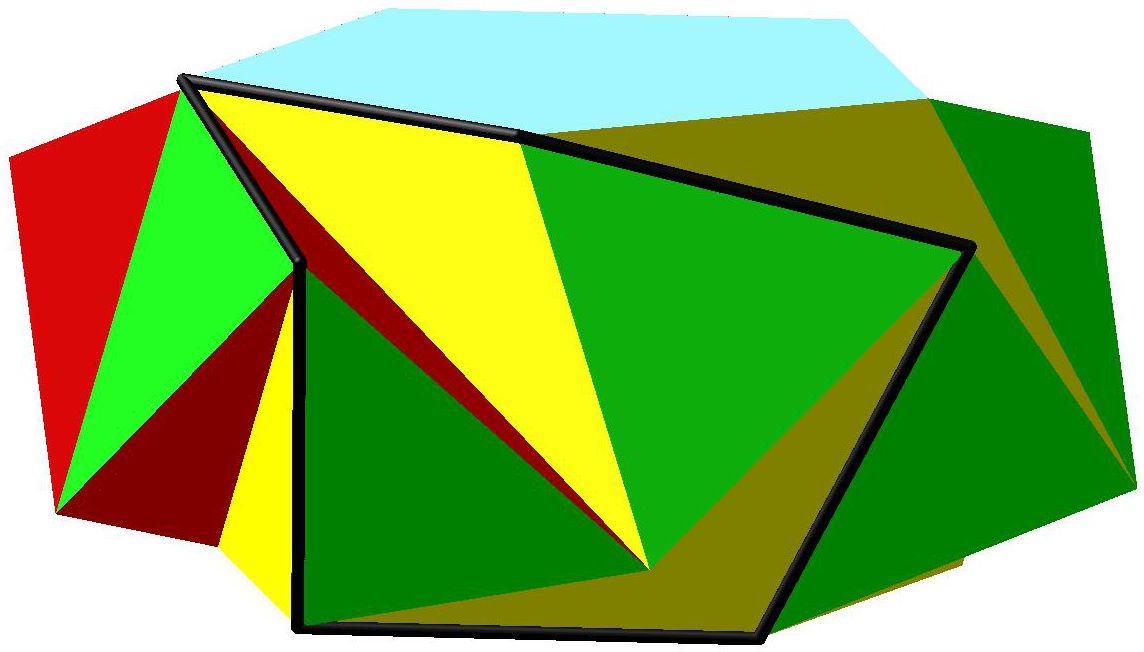}
\begin{small}
\put(0,0){b)}
\end{small}         
  \end{overpic} 
	\hfill
\begin{overpic}
    [width=40mm]{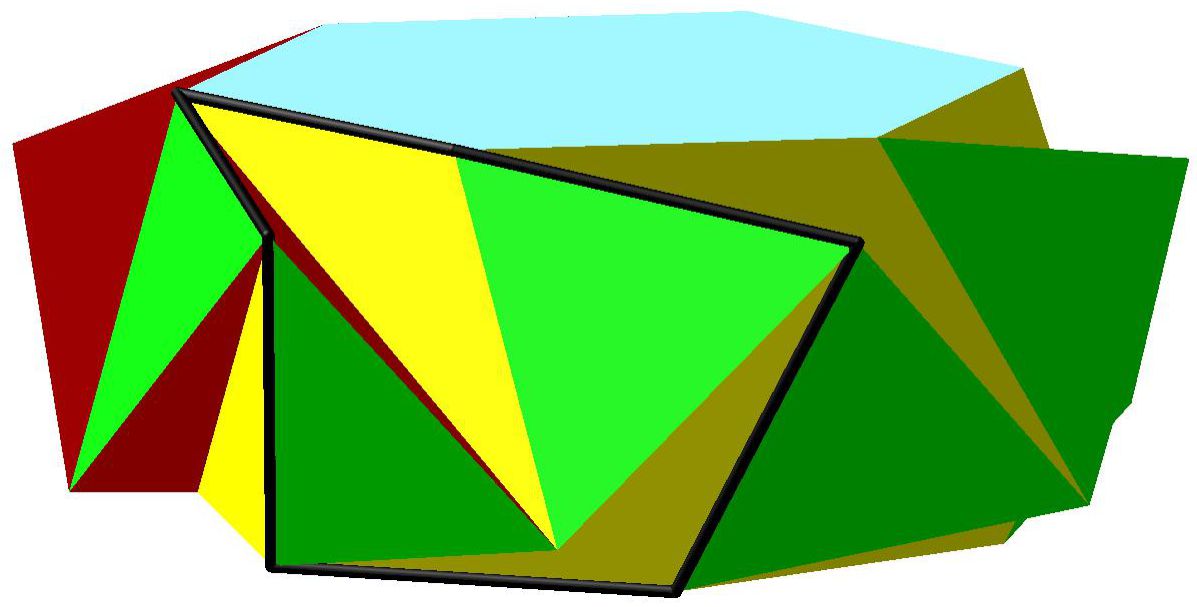}
\begin{small}
\put(0,0){c)}
\end{small}         
  \end{overpic} 	
\caption{Extreme birosettes for $n=4,5,6$ (a,b,c) with yellow/red pedals, where the unit-cells are framed in black.   
}
  \label{fig4}
\end{figure}


\section{Preliminary considerations}\label{sec:prelimi}

Assume that we have given a polyhedron with the combinatorial structure of a birosette for arbitrary $n\geq 3$. 
In order to clarify its degree of freedom (dof), we use the notation of Fig.\ \ref{fig4}a, 
where the vertices of a unit-cell\footnote{A repetitive rotation of the unit-cell around the axis $AB$ about the angle $\tfrac{\pi}{n}$ 
generates the polyhedral belt and increases the indices of the vertices by one (mod $n$).} of the belt are labeled.
The relative position of $\beta$ with respect to $\alpha$ has 6 dofs and each vertex $C_i$ and $D_i$ has further 3 dofs, 
which yields in total $6+6n$ dofs. Moreover, we have $8n$ distance constraints implied by the edges of the structure, with exception of 
those belonging to the skeleton. Therefore the dof of the structure is computed by $F:=6(n+1)-8n$ 
yielding for  $n\geq 3$ a value of $F\leq 0$, which shows that these structures are in general rigid\footnote{Note 
that taking the symmetry of the birosette into account yields $F= 0$ for all $n\geq 3$.}.  

\begin{remark}
It is an open question whether flexible polyhedrons exist having the same combinatorial 
structure as a birosette with $n\geq 3$. 
\hfill $\diamond$
\end{remark}

The assumptions on the edge lengths (inside a unit-cell) of a birosette can be weaken as follows without destroying its symmetry: 
\begin{equation}\label{eq:general}
L_1:=\overline{B_{0}D_{0}}=\overline{A_{1}C_{1}}, \,\,
L_2:=\overline{B_{0}C_{1}}=\overline{A_{0}D_{0}}, \,\,
L_3:=\overline{D_{0}C_{1}}=\overline{C_{1}D_{1}}, \,\,
L_4:=\overline{B_{0}D_{1}}=\overline{A_{0}C_{1}}
\end{equation}
beside the unit length of $\overline{A_{0}A_{1}}$ and $\overline{B_{0}B_{1}}$. Moreover, in the remainder of the 
paper we assume that these skeleton edges are undeformable under the model flexibility in contrast to the 
other edge lengths $L_1,\ldots ,L_4>0$, which can vary\footnote{Note that the study \cite{birosettes} on  
birosettes is more restrictive as only the length $p$ of the diagonal is allowed to change.}.

One can compute (cf.\ Appendix \ref{app:ori} for details) the condition for the developability of the obtained generalized birosette belt in terms of  $L_1,\ldots ,L_4$. Under consideration of the sandglass condition $L_1=L_4$ 
this origami condition reads as: 
\begin{equation}\label{origami}
Q_3=Q_1+Q_2-\sqrt{Q_2(4Q_1-1)} \quad \text{with} \quad Q_i:=L_i^2
\end{equation}

\begin{remark}
Wunderlich further assumed in the study of the trisymmetric sandglass \cite{sanduhr} 
that $L_2=L_3$ holds, which simplifies the related equations considerably and 
allows a very compact treatment. We do not make this assumption here.  \hfill $\diamond$
\end{remark}


\section{Snapping quasi-mechanisms}\label{sec:snap}

As the generalized Wunderlich sandglasses have the same symmetry as the birosettes, the essential vertices of the unit-cell can 
be coordinatized similarly to \cite{birosettes} by
\begin{equation*}
A_0=(R,0,H),  \quad B_0=(Rc,Rs,-H), \quad D_0=(r,0,-h), \quad C_1=(rc,rs,h),
\end{equation*}
using the abbreviations $c:=\cos{\tfrac{\pi}{n}}$ and $s:=\sin{\tfrac{\pi}{n}}$, respectively. 
Moreover, $R$ equals $\tfrac{1}{2s}$ due to the unit-length of the skeleton edges.
With this parametrization we only remain with the following three equations:
\begin{equation}\label{eq:realize}
q_1:\,Q_1-\overline{B_{0}D_{0}}^2=0, \quad
q_2:\,Q_2-\overline{B_{0}C_{1}}^2=0, \quad
q_3:\,Q_3-\overline{D_{0}C_{1}}^2=0, \quad
\end{equation}
as the fourth one, namely $Q_4-\overline{B_{0}D_{1}}^2=0$, is identical with $q_1$. Thus for a given set of squard 
edge lengths $Q_1,Q_2,Q_3$ the corresponding realizations can be computed by solving  $q_1,q_2,q_3$ 
for $H,h,r$, which is straight forward. 

As already mentioned, we want to restrict to structures possessing an {\it extremal snap}, i.e.\ one 
realization is on the boundary of self-intersection. Such configurations are illustrated in Fig.\ \ref{fig6} (right) and have the 
advantage that the self-blocking of the faces increases the structure's load carrying capacity. 
Moreover, the self-covered areas can be provided with holes. In this way the configuration of Fig.\ \ref{fig6} (right)
is still tight, thus it is called {\it closed}, in contrast to the {\it open} one of Fig.\ \ref{fig6} (left). 
Therefore, such structures can for example be used as {\it pressure relief valves}.

\begin{figure}[t]
\hfill
 \begin{overpic}
     [width=25mm]{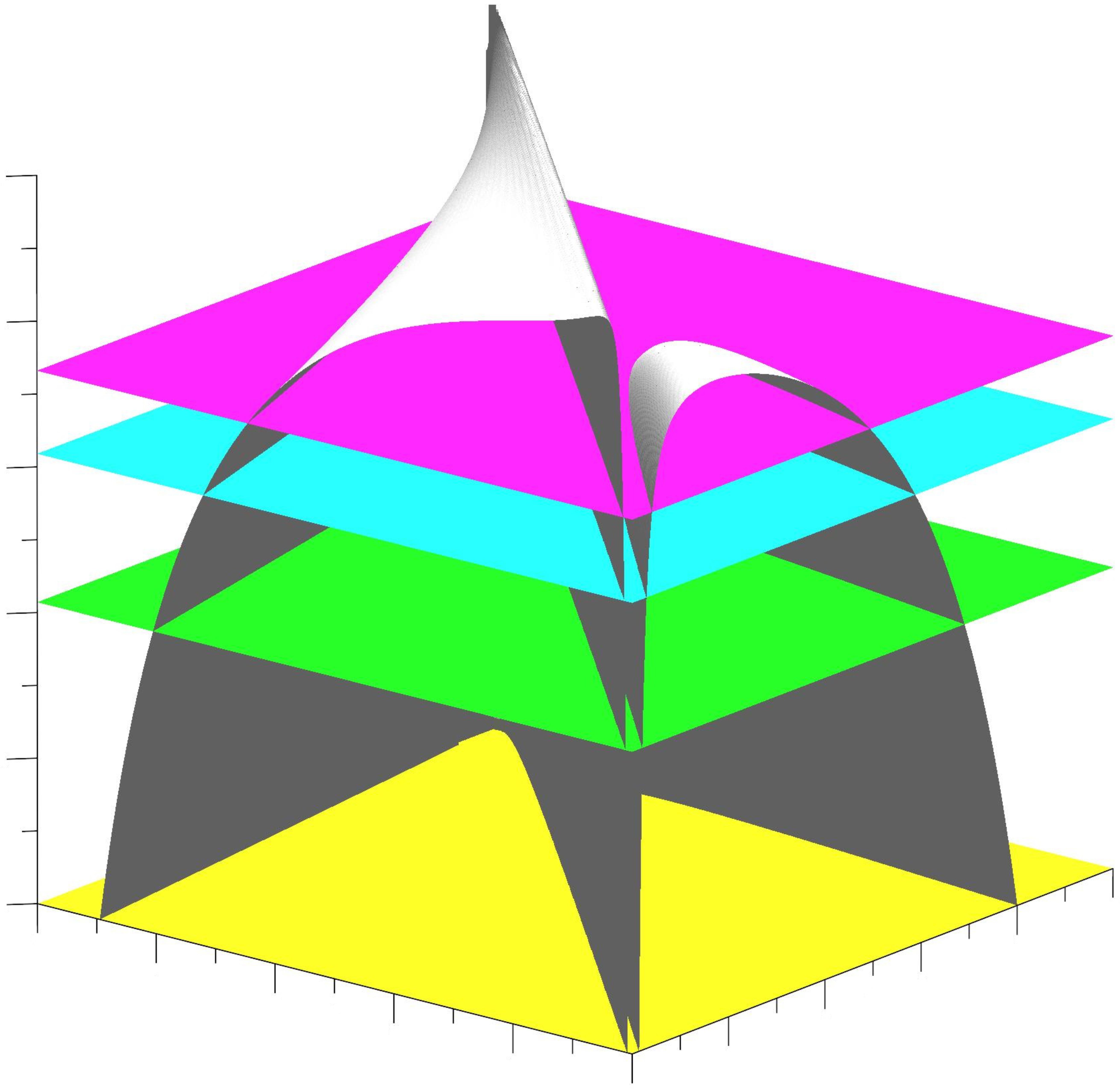}
\begin{small}
\put(0,90){a)}
\put(1,4){$0$}
\put(-7,13){$\tfrac{1}{2}$}
\put(-6,78){$1$}
\put(-6,47){$c$}
\put(11,1){$1$}
\put(23,0){$Q_2$}
\put(73,0){$Q_1$}
\put(97,8){$0$}
\put(88.5,4){$1$}
\end{small}         
  \end{overpic} 
	\hfill
\begin{overpic}
    [width=25mm]{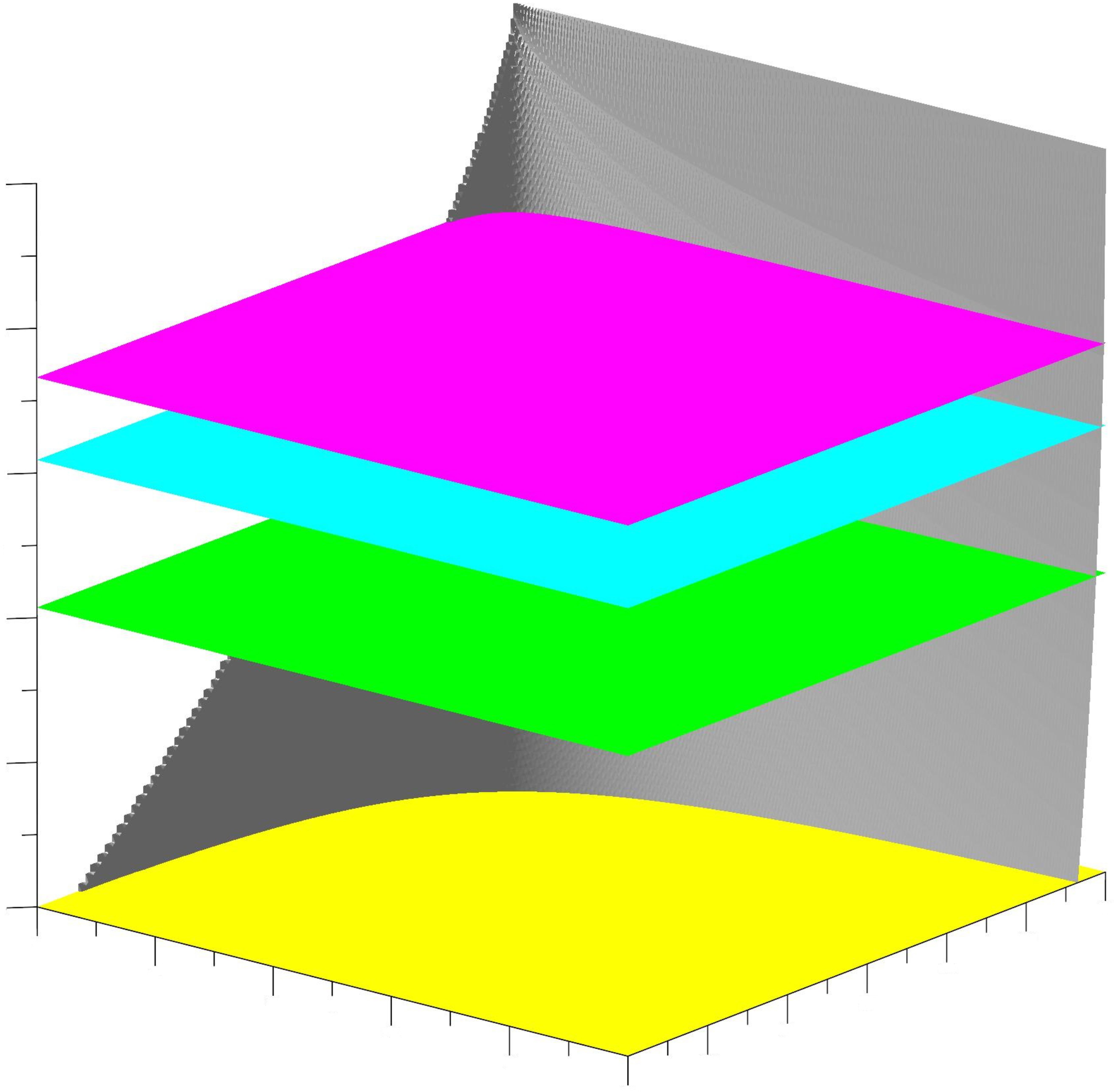}
\begin{small}
\put(0,90){b)}
\put(1,4){$0$}
\put(-7,13){$\tfrac{1}{2}$}
\put(-6,78){$1$}
\put(-6,47){$c$}
\put(11,1){$1$}
\put(23,0){$Q_2$}
\put(71,0){$Q_1$}
\put(97,6){$\tfrac{1}{4}$}
\put(86,3){$\tfrac{26}{10}$}
\end{small}         
  \end{overpic} 	
	\hfill
\begin{overpic}
    [width=25mm]{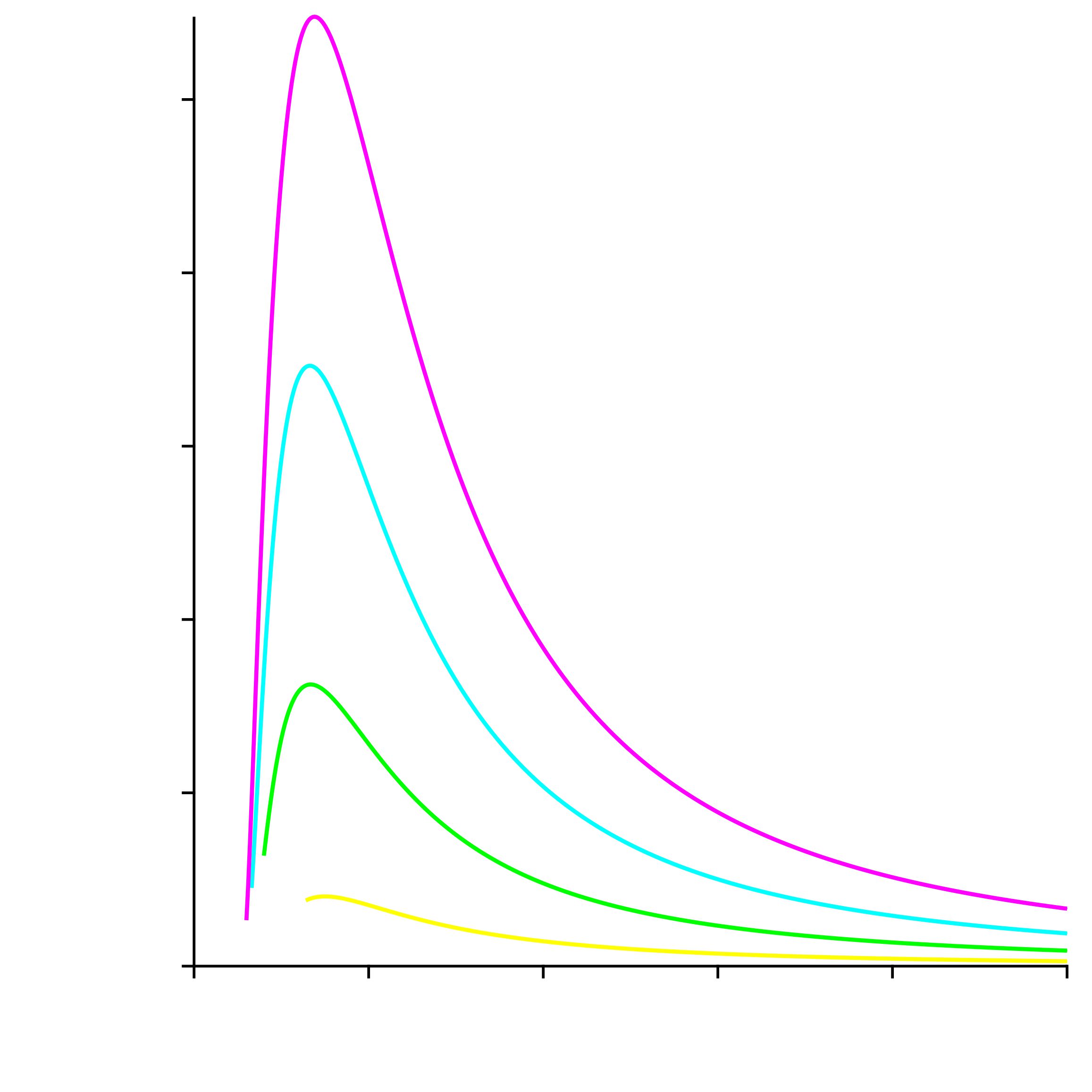}
\begin{small}
\put(0,90){c)}
\put(15,0){$0$}
\put(31,0){$1$}
\put(-9,20){$10^{-4}$}
\put(9,7){$0$}
\put(7,48){$\sigma$}
\put(52,1){$Q_1$}
\end{small}     
  \end{overpic} 
		\hfill
\begin{overpic}
    [width=25mm]{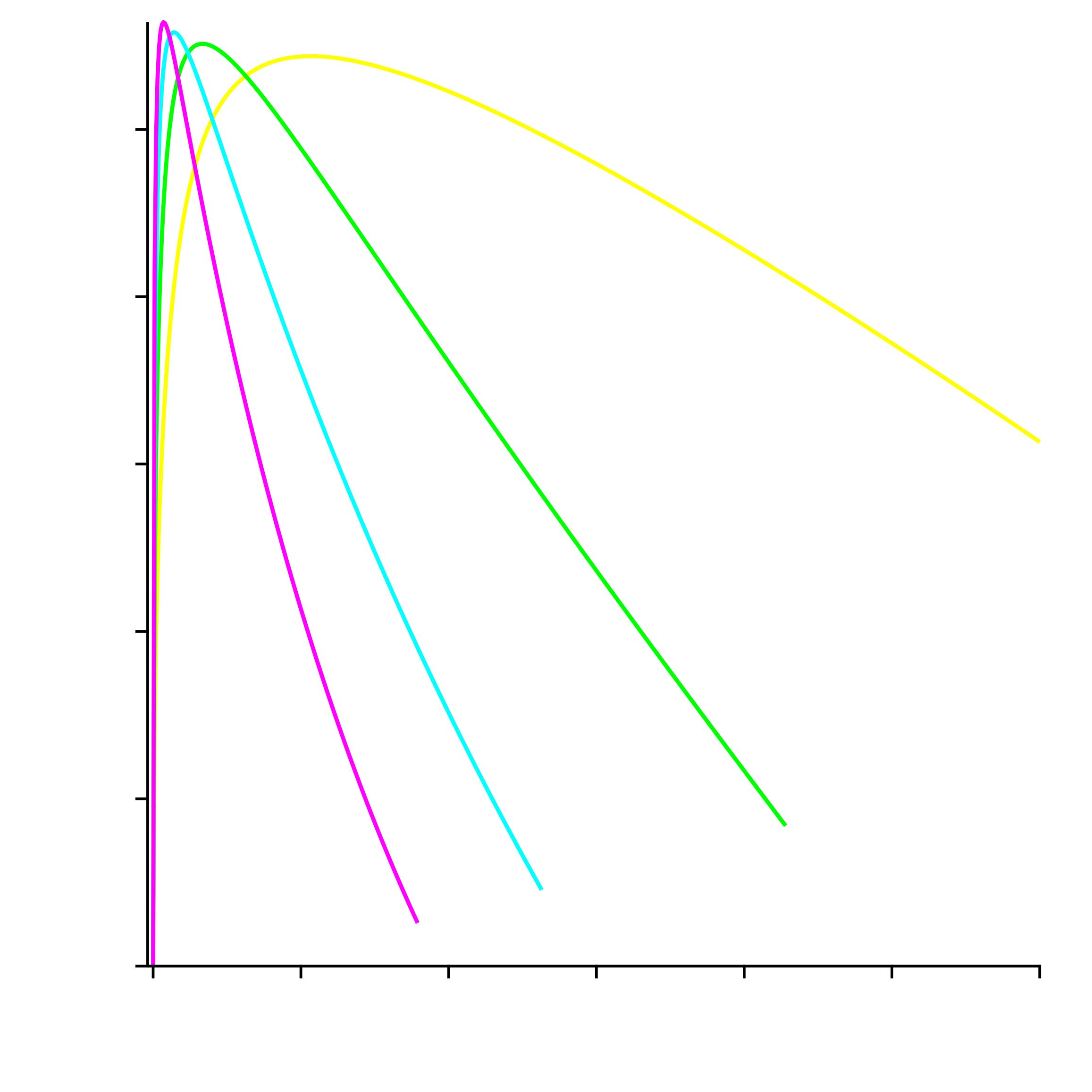}
\begin{small}
\put(0,90){d)}
\put(11,-1){$\tfrac{1}{4}$}
\put(21.5,-1){$\tfrac{26}{10}$}
\put(-10,22){$0.05$}
\put(5,7){$0$}
\put(4,47){$\kappa$}
\put(58,0){$Q_1$}
\end{small}     
  \end{overpic}
\caption{
Visualization of (a) Eq.\ (\ref{eq:extrem}) and (b) the shakiness condition, which only yields feasible values for $Q_2$ within 
a narrow domain of $Q_1$, where $Q_1=\tfrac{1}{4}$ is an asymptote for the $Q_2$ values.  The horizontal $c$-planes for 
$n=3,\ldots,6$ are colored in yellow, green, cyan and magenta, respectively. The same color-coding is used for the graphs of the snappability (c) and 
the shakeability (d). 
}
  \label{fig5}
\end{figure}

\begin{figure}[t]

\hfill
 \begin{overpic}
    [height=16mm]{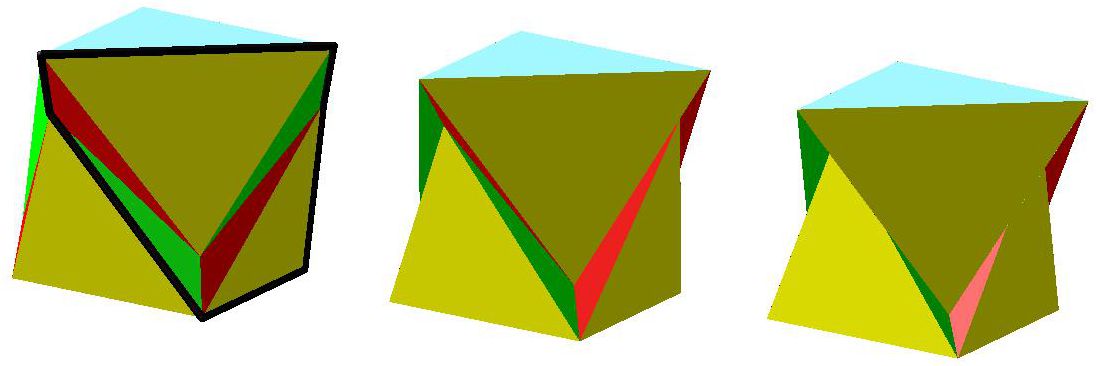}
\begin{small}
\put(0,0){a)}
\end{small}         
  \end{overpic} 
	\hfill
\begin{overpic}
    [height=16mm]{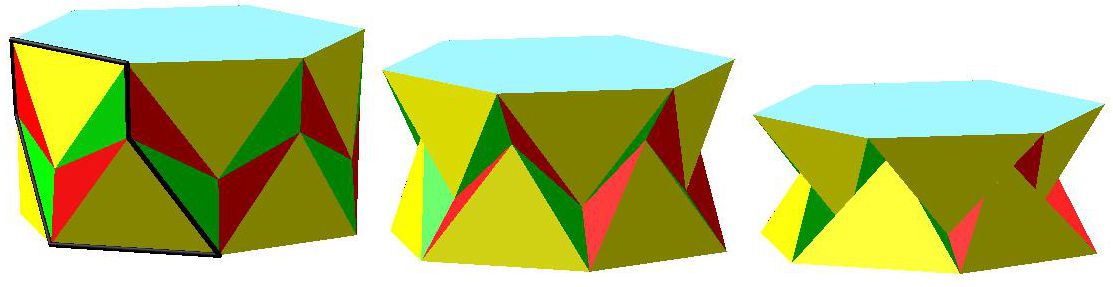}
\begin{small}
\put(-1,0){d)}
\end{small}         
  \end{overpic} 
	\hfill \smallskip
	\newline 
 \begin{overpic}
    [height=16mm]{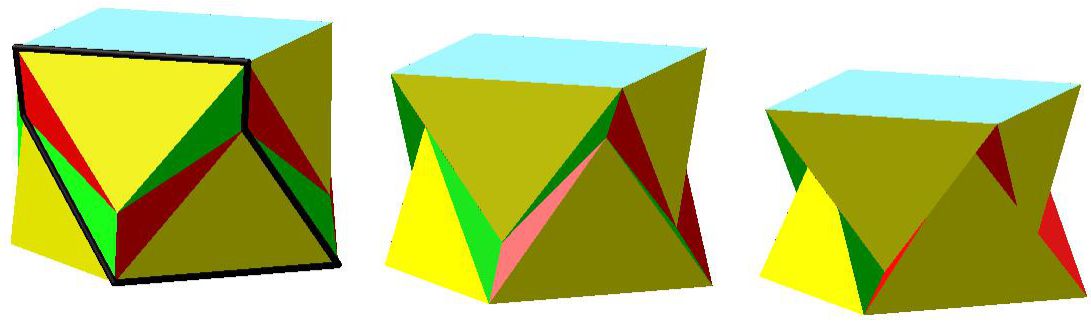}
\begin{small}
\put(0,0){b)}
\end{small}         
  \end{overpic} 
	\hfill
\begin{overpic}
    [height=16mm]{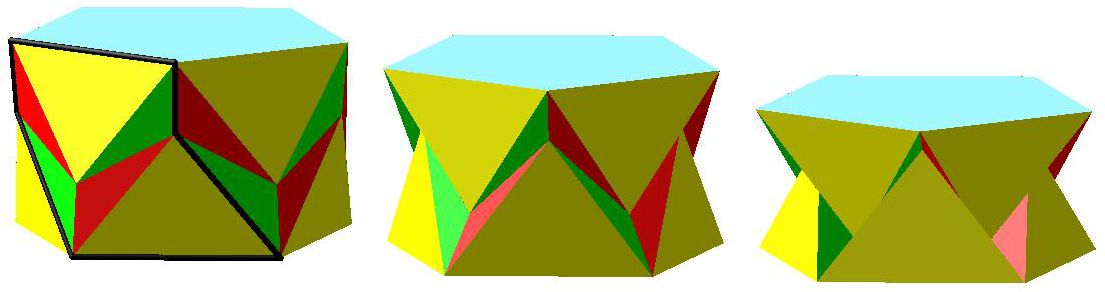}
\begin{small}
\put(0,0){c)}
\end{small}         
  \end{overpic} 
\caption{Illustration of the snap from the open state (left) over the passed shaky configuration (center) to the closed state (right) 
for the structures causing the maximal snappability for $n=3$ (a), $n=4$ (b), $n=5$ (c) and $n=6$ (d). Moreover, the unit-cells of the 
open state are framed in black. 
}
  \label{fig6}
\end{figure}

At the closed state the points $A_0,B_0,D_0,C_1$ are coplanar. This is the case for $(2rs - 1)(2Hrs + h)=0$ where the factors imply 
a dihedral angle of $\pi$ and $0$, respectively, along the edge $D_0C_1$. 
Therefore we set  $r=-\tfrac{h}{2Hs}$ and plug this expression into $q_1,q_2,q_3$. 
By eliminating the unknowns $H,r,s$ from these three equations together with $c^2+s^2-1=0$  by means of resultants, we end up with
\begin{equation} \label{eq:extrem}
\begin{split}
&4cQ_2Q_1 - 2cQ_2^2 - 2Q_1^2 - 28Q_2Q_1 - 2Q_2^2 + Q_1 + 5Q_2- 2cQ_1^2 \\
&+W^{3/2}\sqrt{Q_2} + 8Q_2^{3/2}\sqrt{W} + 4Q_1\sqrt{Q_2}\sqrt{W}=0 
\end{split}
\end{equation}
by taking into account Eq.\ (\ref{origami}) and $W:=4Q_1 - 1$. This equation is linear\footnote{For $Q_1=Q_2$ it is independent of $c$ but then it can only vanish for $Q_1=Q_2=0$, a contradiction. } in $c$ and plotted in Fig.\ \ref{fig5}a. 
Moreover, it can be solved explicitly for $Q_2$ (e.g.\ with {\sc Maple}), which yields four branches. We let $Q_1$ run 
within the interval $]0.25;5]$ in steps of $0.01$ for $n=3,\ldots ,6$ and check if this value implies a pair of snapping sandglass realizations. 
To do so, we compute the realizations and check in a first step if they are free of self-intersections, and in a second step 
if they can snap into each other by the methods presented by the author in \cite{nawr1,nawr2}. 
It turns out that only the branch of Eq.\ (\ref{eq:extrem}), which yields the lowest value for $Q_2$ carries solutions of our problem. 
Animations of these families of snapping sandglasses (for $n=3,\ldots ,6$) can be downloaded from the 
author's \medskip homepage\footnote{\label{homepage}{\url{https://www.dmg.tuwien.ac.at/nawratil/publications.html}}}.

\noindent
{\bf Snappability.} According to \cite{nawr1,nawr2} the snap between two sandglass realizations has to pass a shaky configuration 
(with squared edge lengths $S_1,S_2,S_3$) at the maximum state of deformation (see Fig.\ \ref{fig6} (center)), which is used for the evaluation of the snapping capability 
in terms of the {\it snappability}. 
By considering the belt as a joint-bar structure\footnote{It can also be interpreted as a panel-hinge framework but the formula is more complicated \cite{nawr1,nawr2}.}
this index $\sigma$, which is based on the total elastic strain energy density of the framework, can be computed as follows (cf.\ \cite{nawr1,nawr2}):
\begin{equation}\label{eq:snap}
\sigma:=\left(4n\tfrac{(Q_1-S_1)^2}{8L_1^3} + 2n\tfrac{(Q_2-S_2)^2}{8L_2^3} + 2n\tfrac{(Q_3-S_3)^2}{8L_3^3}\right)/(4nL_1+2nL_2+2nL_3)
\end{equation}
In Fig.\ \ref{fig5}c the snappability of the computed snapping sandglass realizations, which are free of self-intersections, are displayed. 
The structures causing the maximal snappability are visualized in Fig.\ \ref{fig6}. 
Additional graphs concerning the change in volume, height or waist radius during the snap are given in Appendix \ref{app:add}
as well as the crease pattern of the structure displayed in Fig.\ \ref{fig6}a.

\section{Shaky quasi-mechanisms}\label{sec:shake}

For the computation of shaky origami-like sandglasses we proceed as follows. We eliminate from $q_1,q_2,q_3$ of 
Eq.\ (\ref{eq:realize}) the unknowns $H,h$ by means of resultant, such that we end up with a polynomial in $r$. 
In order that this polynomial has a solution of higher order, its discriminant with respect to $r$ has to vanish.
From the resulting expression we eliminate $s$ by applying the resultant method with respect to $c^2+s^2-1=0$. 
Under consideration of Eq.\ (\ref{origami}) we end up with the shakiness condition of the form $w_4c^4+w_3c^3+w_2c^2+w_1c+w_0=0$ 
whose coefficients are given in Appendix \ref{app:coef} due to their length.  
This quartic equation in $c$ is plotted in  Fig.\ \ref{fig5}b.  
Note that this equation cannot be solved explicitly for $Q_2$ but we can evaluate it \medskip numerically.

The  shakiness of the structure is in general related with a non-trivial infinitesimal isometric deformation of dimension 1. 
The corresponding velocity vectors $\Vkt v(X_i)$ of the vertices $X_i$  ($i=0,\ldots ,n-1$) with $X\in\left\{A,B,C,D\right\}$ 
are determined\footnote{\label{note:inf}Note that the faces $\alpha$ and $\beta$ are assumed to 
translate instantaneously along the $z$-axis with the same speed but in opposite direction; i.e.\ $\Vkt v(A_i)=-\Vkt v(B_i)$. 
Therefore the symmetry of the structure remains intact by attaching the velocity vectors $\Vkt v(X_i)$ to the corresponding vertices $X_i$.} 
up to a non-zero factor (fixed by Eq.\ (\ref{normalize}) given later on) and can be computed by applying the 
{\it projection theorem}. For details the interested reader is referred to \medskip  Appendix \ref{app:inf}.

\begin{figure}[t]
\hfill
\begin{overpic}
    [width=20mm]{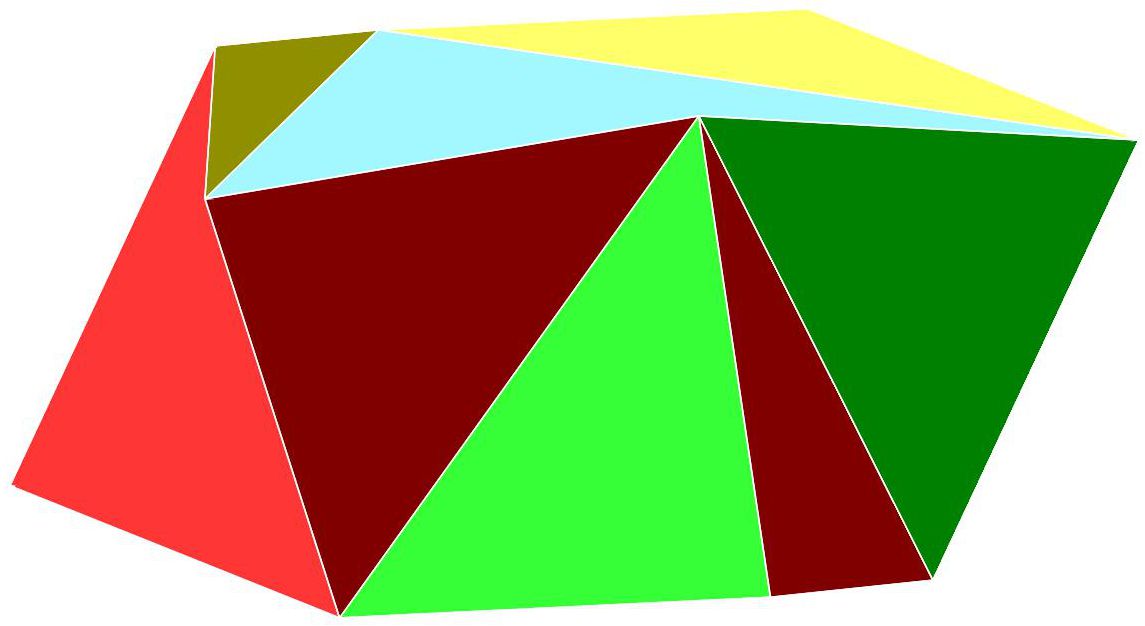}
\begin{small}
\put(-7,0){a)}
\end{small}     
  \end{overpic} 
\hfill
 \begin{overpic}
    [width=23mm]{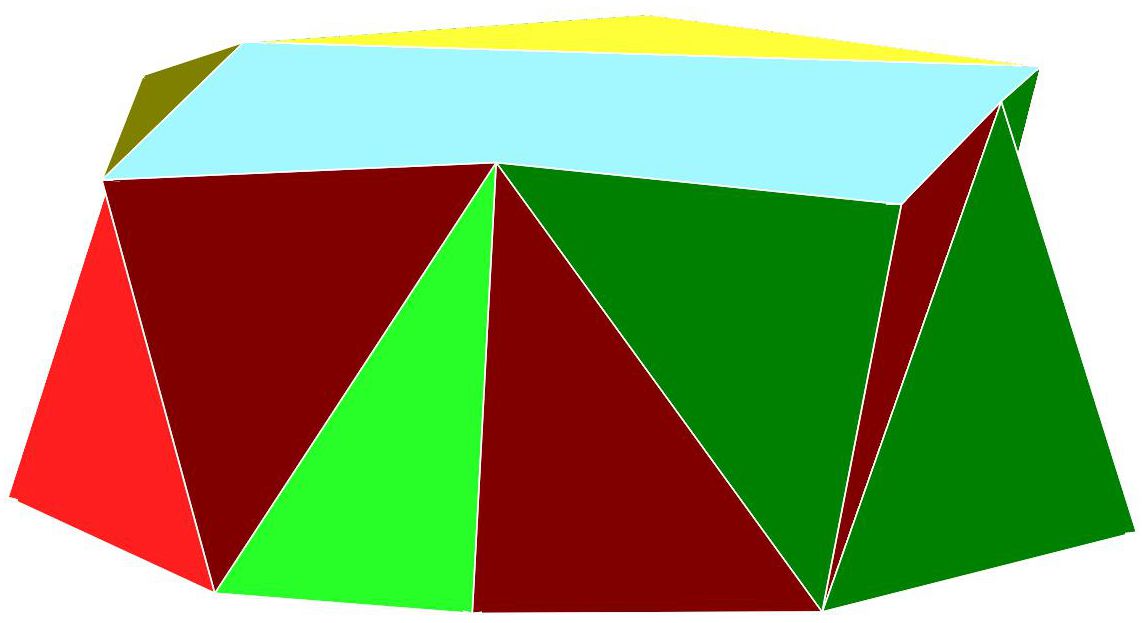}
\begin{small}
\put(-6,0){b)}
\end{small}         
  \end{overpic} 
	\hfill
\begin{overpic}
    [width=30mm]{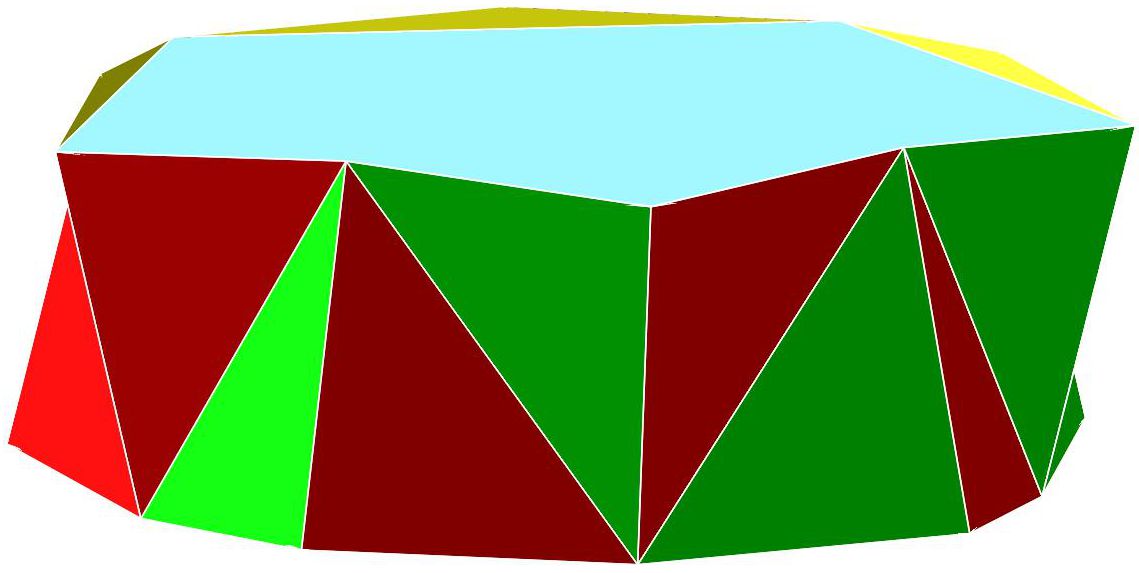}
\begin{small}
\put(0,0){c)}
\end{small}         
  \end{overpic} 	
\hfill
\begin{overpic}
    [width=35mm]{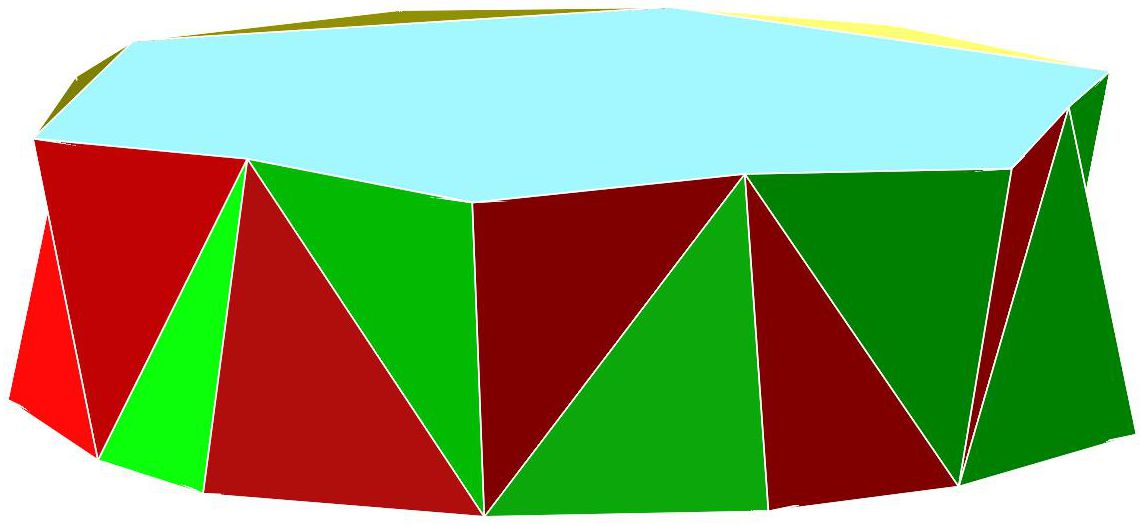}
\begin{small}
\put(0,0){d)}
\end{small}         
  \end{overpic} 		
\caption{
Shaky structures causing the maximal shakeability for $n=3,4,5,6$ (a,b,c,d).
}
  \label{fig7}
\end{figure}

\noindent
{\bf Shakeability.} In order to evaluate the structure's capability to shake, we introduce the so-called 
{\it shakeability $\kappa$}, which is defined as curvature of the snappability function over the space of squared edge lengths in direction associated with the infinitesimal mobility.
As the snappability function is already dimensionless we also have to normalize the velocity vectors in such a way.
This can e.g.\ be achieved by the condition\footnote{\label{foot:con}One can also think of other normalizations, which will effect the resulting value for $\kappa$.}, 
that the mean of the relative instantaneous changes of the squared edge 
lengths is equal to 1. For our sandglass structure this normalization reads as follows:
\begin{equation}\label{normalize}
\left(4n\tfrac{\|\Vkt v(B_0)-\Vkt v(D_0)\|^2}{Q_1} +  2n\tfrac{\|\Vkt v(B_0)-\Vkt v(C_1)\|^2}{Q_2} + 2n\tfrac{\|\Vkt v(D_0)-\Vkt v(C_1)\|^2}{Q_3}\right)/(8n) =1
\end{equation}
Assumed that this condition holds, we can set 
$S_1=Q_1+t\|\Vkt v(B_0)-\Vkt v(D_0)\|^2$, 
$S_2=Q_2+t\|\Vkt v(B_0)-\Vkt v(C_1)\|^2$ and 
$S_3=Q_3+t\|\Vkt v(D_0)-\Vkt v(C_1)\|^2$
and plug these expressions into Eq.\ (\ref{eq:snap}), which now depends quadratically on $t$; i.e.\ $\sigma(t)$.  
According to the well-known curvature formula the shakeability $\kappa$ can then be computed as
\begin{equation}
\kappa:=\left.\tfrac{\sigma''}{(1+\sigma'^2)^{3/2}}\right|_{t=0}=\left. \sigma''\right|_{t=0}
\end{equation}
We let $Q_1$ run within the interval $]0.25;0.31]$ in steps of $0.001$ for $n=3,\ldots ,6$ and 
compute for the associated shaky realization the shakeability, which is displayed in Fig.\ \ref{fig5}d.  
The structures causing the maximal shakeability are visualized in Fig.\ \ref{fig7}. 

Animations of the families of shaky origami-like sandglasses for $n=3,\ldots ,6$ 
can be downloaded from the author's homepage (cf.\ footnote \ref{homepage}) 
and the crease pattern of the structure displayed in Fig.\ \ref{fig7}a is given in Appendix \ref{app:add}. 

\begin{remark}
Shaky realizations with smaller $\kappa$ are more shaky, but note that 
the shakeability cannot vanish as in this case all the velocity vectors have to be the same; i.e.\ 
the infinitesimal isometric deformation is trivial (instantaneous translation).  \hfill $\diamond$
\end{remark}

\section{Conclusion and future work}\label{sec:end}

We generalized Wunderlich's sandglass polyhedron in analogy to the birosette construction, with the
additional feature that the belt of the antiprismatic skeleton is developable.  One can think of two further generalizations 
of the obtained snapping or shaky quasi-mechanisms; namely 
(a) to omit the sandglass condition and/or the origami condition for the study of the generalized birosette structures (given by Eq.\ (\ref{eq:general}))
and (b) to use an antifrustum as skeleton (i.e.\ $\alpha$ and $\beta$ have different radii). 
For the latter case one can proceed similar to \cite{conformal,nawr3}. 
Moreover, we introduced the {\it shakeability} to evaluate the capability of a shaky quasi-mechanism to flex 
as counterpart to the snappability \cite{nawr1,nawr2}. 
As this index is dimensionless, it enables a comparison of shaky structures differing in the inner geometry, 
which is also subject to future research as well as the influence of the normalization condition (cf.\ footnote \ref{foot:con}).  

\begin{acknowledgement}
The author is supported by grant P\,30855-N32 of the Austrian Science Fund FWF 
and by FWF project F77 (SFB ``Advanced Computational Design'', subproject SP7).
\end{acknowledgement}

\section*{Appendix}

\subsection{Origami condition for the generalized birosette belt}\label{app:ori}

\begin{figure}[b]
$\phm$
\hfill
\begin{overpic}
    [width=90mm]{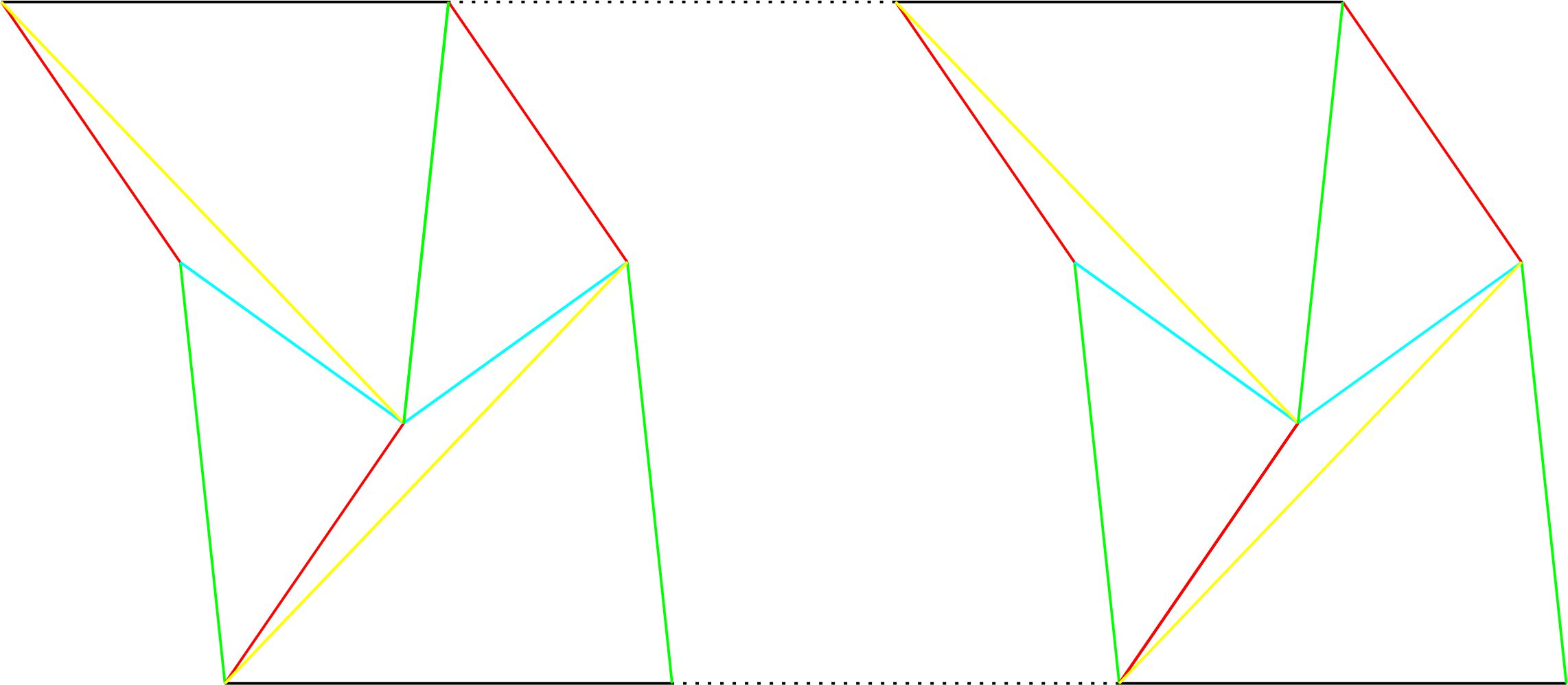}
\begin{small}
\put(10,1.5){$B_0^*$}
\put(43.5,1.5){$B_1^*$}
\put(64,1.5){$B_{n-1}^*$}
\put(101,1.5){$B_0^*$}
\put(-3,40){$A_0^*$}
\put(24,40){$A_1^*$}
\put(52,40){$A_{n-1}^*$}
\put(81,40){$A_0^*$}
\put(7,25){$D_0^*$}
\put(41,25){$D_1^*$}
\put(61,25){$D_{n-1}^*$}
\put(98,25){$D_0^*$}
\put(78,15.5){$C_0^*$}
\put(20.5,15.5){$C_1^*$}
\end{small}     
  \end{overpic} 
\hfill $\phm$
\caption{Sketch for the computation of the origami condition for the generalized birosette belt. 
Edges with the same length are colored equally (cf.\ Eq.\ (\ref{eq:dev})).
}
  \label{fig:net}
\end{figure}
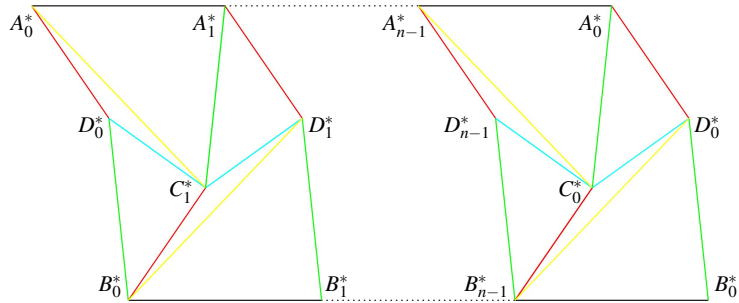

In this section we compute the condition for the developability of the obtained generalized birosette belt in terms of  $L_1,\ldots ,L_4$. 
We start with a developed unit-cell, whose  vertices can be coordinatized with respect to a planar Cartesian frame as:
\begin{equation}
\begin{split}
&B_0^*=(0,0), \quad B_1^*=(1,0), \quad D_0^*=(a,b), \quad D_1^*=(a+1,b),  \\
&C_1^*=(c,d), \quad A_0^*=(a+e,b+f), \quad A_1^*=(1+a+e+1,b+f),
\end{split}
\end{equation}
without loss of generality.  Note that the repetition of the unit-cell implies the parallelism of $B_0^*D_0^*\parallel B_1^*D_1^*$ 
as well as $A_0^*D_0^*\parallel A_1^*D_1^*$ (cf.\ Fig.\ \ref{fig:net}), where the upper index $*$ indicates the development. 
Moreover, we have to assume that $b,d,f>0$ has to hold in order to avoid overlaps of the development.  
It can easily be verified that the conditions (cf.\ Eq. (\ref{eq:general}))
\begin{equation}\label{eq:dev}
\overline{B_{0}^*D_{0}^*}=\overline{A_{1}^*C_{1}^*}, \quad
\overline{B_{0}^*C_{1}^*}=\overline{A_{0}^*D_{0}^*}, \quad
\overline{D_{0}^*C_{1}^*}=\overline{C_{1}^*D_{1}^*}, \quad
\overline{B_{0}^*D_{1}^*}=\overline{A_{0}^*C_{1}^*}
\end{equation}
imply the relations $d=f$, $e=-c$ and  $c=a+\tfrac{1}{2}$.
Now the remaining unknowns $a,b,f$ can be expressed in terms of $L_1,L_2,L_4$ from 
\begin{equation}
Q_1-\overline{B_{0}^*D_{0}^*}^2=0, \quad
Q_2-\overline{B_{0}^*C_{1}^*}^2=0, \quad
Q_4-\overline{B_{0}^*D_{1}^*}^2=0
\end{equation}
with $Q_i:=L_i^2$, which yields
\begin{equation*}
a = \tfrac{Q_4-Q_1-1}{2},\quad   
b = \tfrac{1}{2}\sqrt{ 2Q_1 + 2Q_4 - 1 -(Q_1-Q_4)^2},\quad 
f = \tfrac{1}{2}\sqrt{4Q_2 -(Q_1-Q_4)^2}
\end{equation*}
Plugging the obtained expressions into $Q_3-\overline{D_{0}^*C_{1}^*}^2=0$ implies the 
origami condition in terms of squared edge lengths. This condition simplifies to Eq.\ (\ref{origami})
if one takes the sandglass condition $L_1=L_4$ ($\Rightarrow$ $Q_1=Q_4$) into account.

\newpage

\subsection{Coefficients of the shakiness condition} \label{app:coef}

The coefficients of the shakiness condition computed in Section  \ref{sec:shake}
read as follows:
\begin{equation}
\begin{split}
w_4=&128Q_1^2Q_2^2 \\
w_3=&\Big[(64-96Q_1)Q_2^{5/2} - 96Q_1^2Q_2^{3/2}\Big]\sqrt{W} + 64Q_1(2Q_1 - 1)Q_2^2 \\
w_2=&96\Big[Q_1^3 + (2Q_2 + \tfrac{5}{16})Q_1^2 + (Q_2^2 - \tfrac{23Q_2}{8})Q_1 - \tfrac{13Q_2^2}{48} + \tfrac{Q_2}{2}\Big]Q_2- \\
&48\sqrt{W}Q_1(2Q_1 - 1)Q_2^{3/2} - 24W^{3/2}Q_2^{5/2} \\
w_1=&96\Big[Q_1^3 + (2Q_2 - \tfrac{3}{4})Q_1^2 + (Q_2^2 - Q_2 - \tfrac{5}{32})Q_1 - \tfrac{Q_2^2}{4} + \tfrac{7Q_2}{32}\Big]Q_2- \Big[ 8Q_2^{7/2}- \\
&(15 + 36Q_1 -24Q_1^2)Q_2^{3/2} + (24Q_1 + 38)Q_2^{5/2}  - (2Q_1^2-8Q_1^3)\sqrt{Q_2}\Big]\sqrt{W} \\
w_0=&  (148Q_2^2 + 35Q_2)Q_1 + 2Q_2^3 - 37Q_2^2 - Q_2 - 70Q_1^2Q_2-\Big[ 8Q_2^{7/2}+ \\
&(24Q_1^2 - 12Q_1 + 24)Q_2^{3/2} + 6WQ_2^{5/2}  + (8Q_1^3 - 6Q_1^2 + Q_1)\sqrt{Q_2}\Big]\sqrt{W} 
\end{split}
\end{equation}

\subsection{Computing the infinitesimal motion} \label{app:inf}

From item (3) of the review part and under consideration of footnote \ref{note:inf},  
the velocity vectors of $A_i$ and $B_i$ can be given as: 
\begin{equation}
\Vkt v(A_i)=(0,0,z), \quad \Vkt v(B_i)=(0,0,-z). 
\end{equation}
Moreover, due to the symmetry of the structure we can set $\Vkt v(D_0)=(u,0,v)$ ($\Rightarrow$  $\Vkt v(C_1)=(uc,us,-v)$).  
These unknowns $u,v$ can be computed in dependency of $z$ from the following two equations, which arise from applying twice the {\it projection theorem}; i.e.\
\begin{equation}
(D_0-A_0)\Vkt v(A_0)=(D_0-A_0)\Vkt v(D_0), \quad 
(D_0-B_0)\Vkt v(B_0)=(D_0-B_0)\Vkt v(D_0).
\end{equation}
This already determines the infinitesimal motion up to the non-zero factor $z$.

\subsection{Additional information}\label{app:add}

\begin{figure}[h]
$\phm$
\hfill
\begin{overpic}
    [width=94mm]{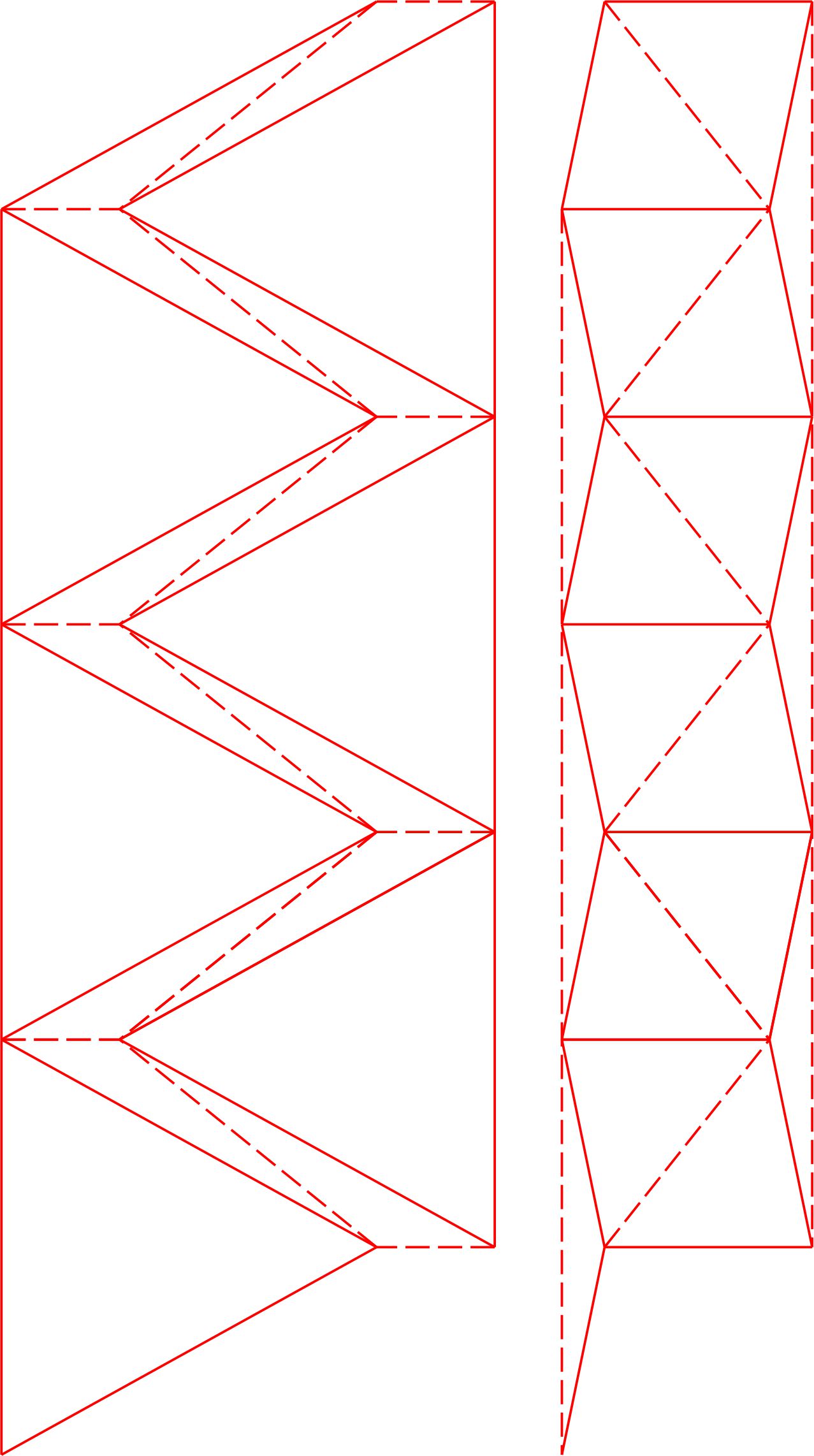}
\begin{small}
\put(4,0){a)}
\put(40,0){b)}
\end{small}     
  \end{overpic} 
\hfill $\phm$
\caption{(a) Crease pattern of the snapping structure displayed in Fig.\ \ref{fig6}a, where  the
mountain folds are solid and the valley folds are dashed with respect to the closed state. 
(b) Crease pattern of the shaky structure displayed in Fig.\ \ref{fig7}a, where the folds are marked in the same way.
}
  \label{fig9}
\end{figure}

\noindent
{\bf Crease patterns.} 
The crease patterns of the structures displayed in Figs.\ \ref{fig6}a  and \ref{fig7}a are given in 
Fig.\ \ref{fig9}. Note that for $n=3$ the  materialization of the faces
$\alpha$ and $\beta$ can be omitted as they are  triangles. 

\newpage

\newpage

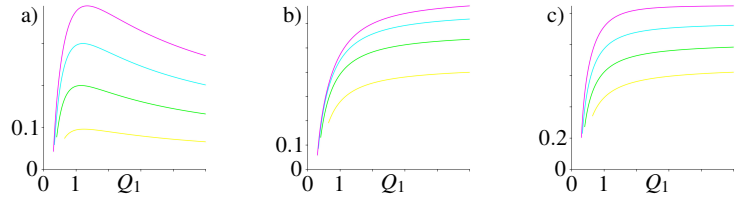
\begin{figure}[h]
\hfill
\begin{overpic}
    [width=25mm]{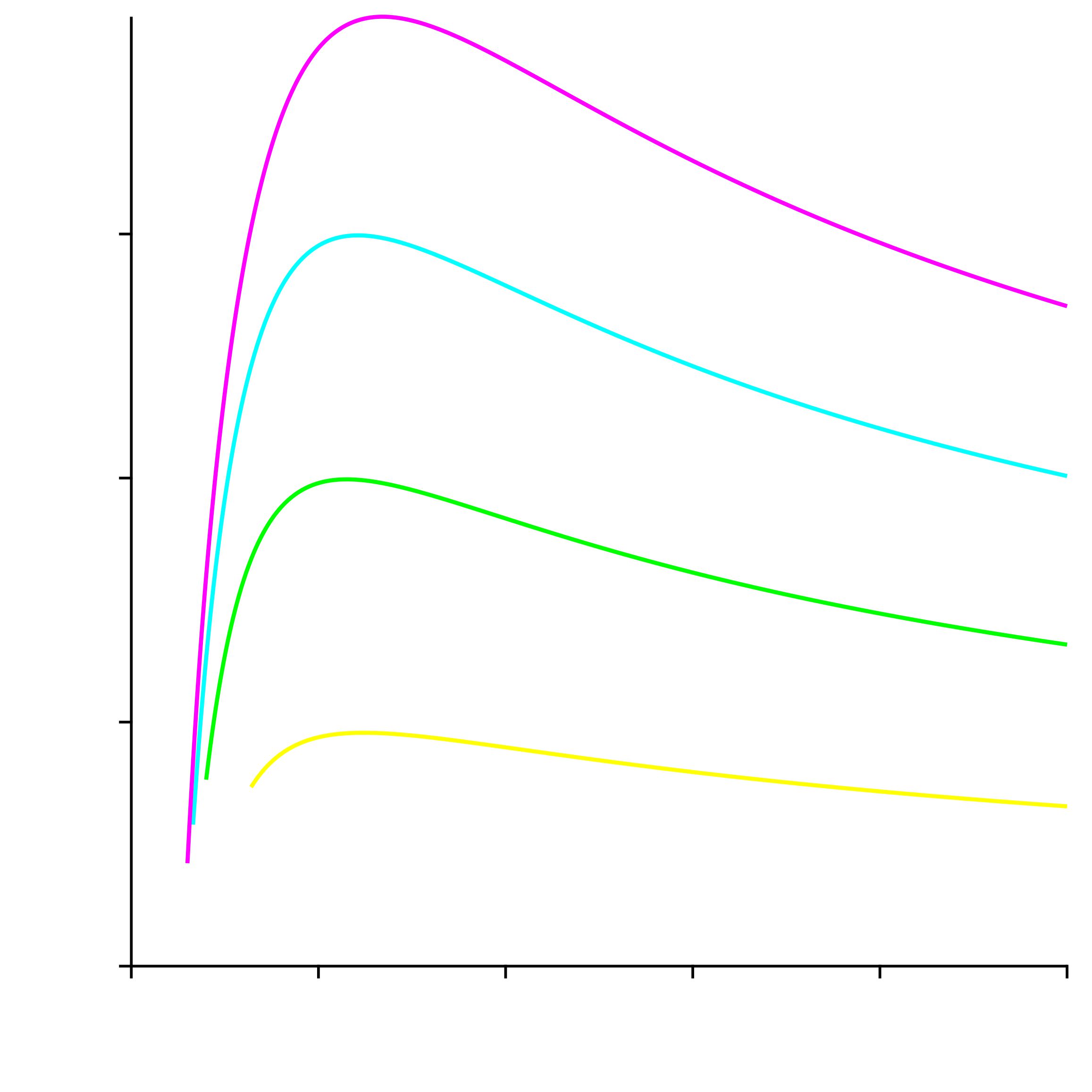}
\begin{small}
\put(0,90){a)}
\put(9,0){$0$}
\put(26,0){$1$}
\put(-5,29.5){$0.1$}
\put(3,8){$0$}
\put(51,1){$Q_1$}
\end{small}         
  \end{overpic} 	
	\hfill
\begin{overpic}
    [width=25mm]{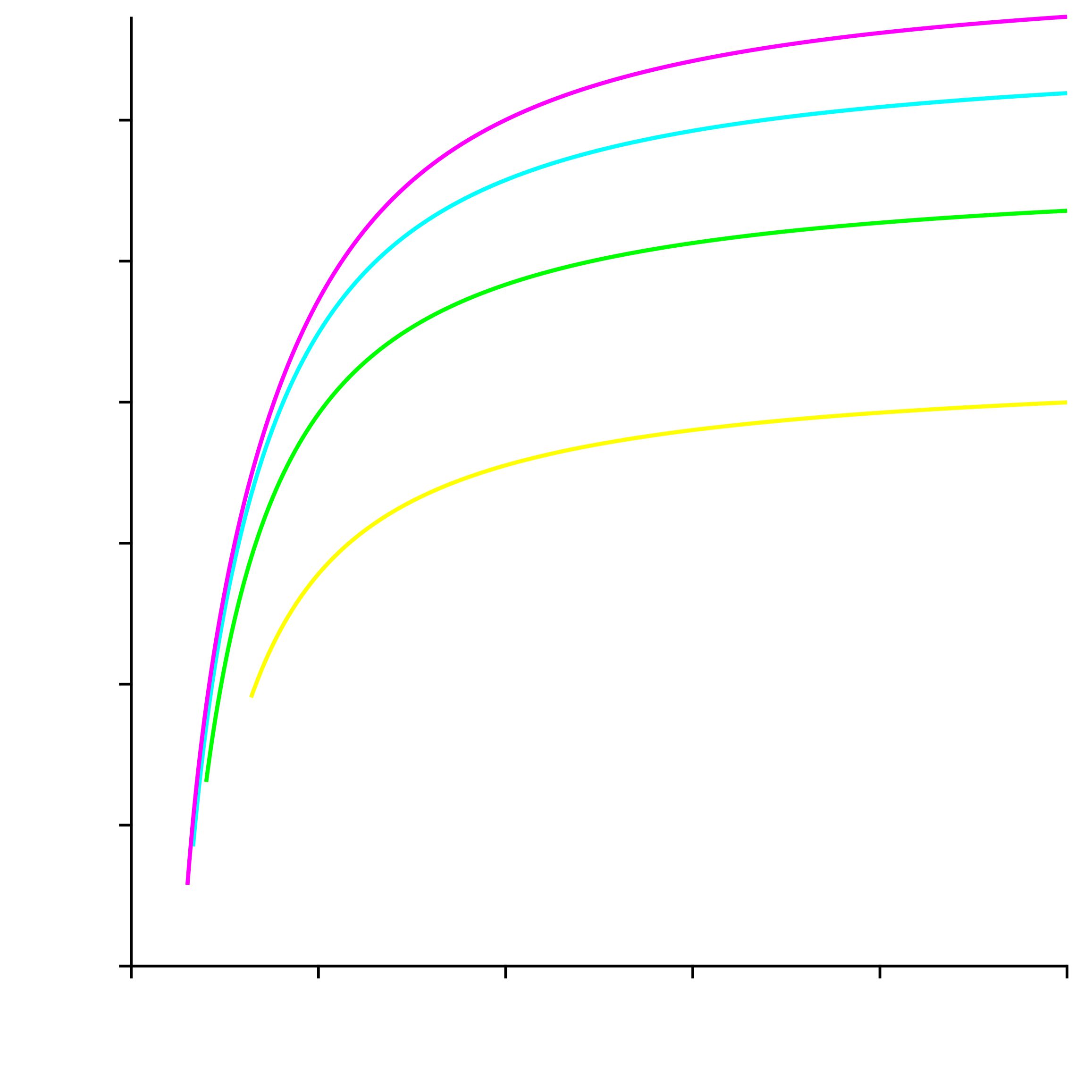}
\begin{small}
\put(-1,90){b)}
\put(9,0){$0$}
\put(26,0){$1$}
\put(-5,20.5){$0.1$}
\put(3,8){$0$}
\put(51,1){$Q_1$}
\end{small}     
  \end{overpic} 
		\hfill
\begin{overpic}
    [width=25mm]{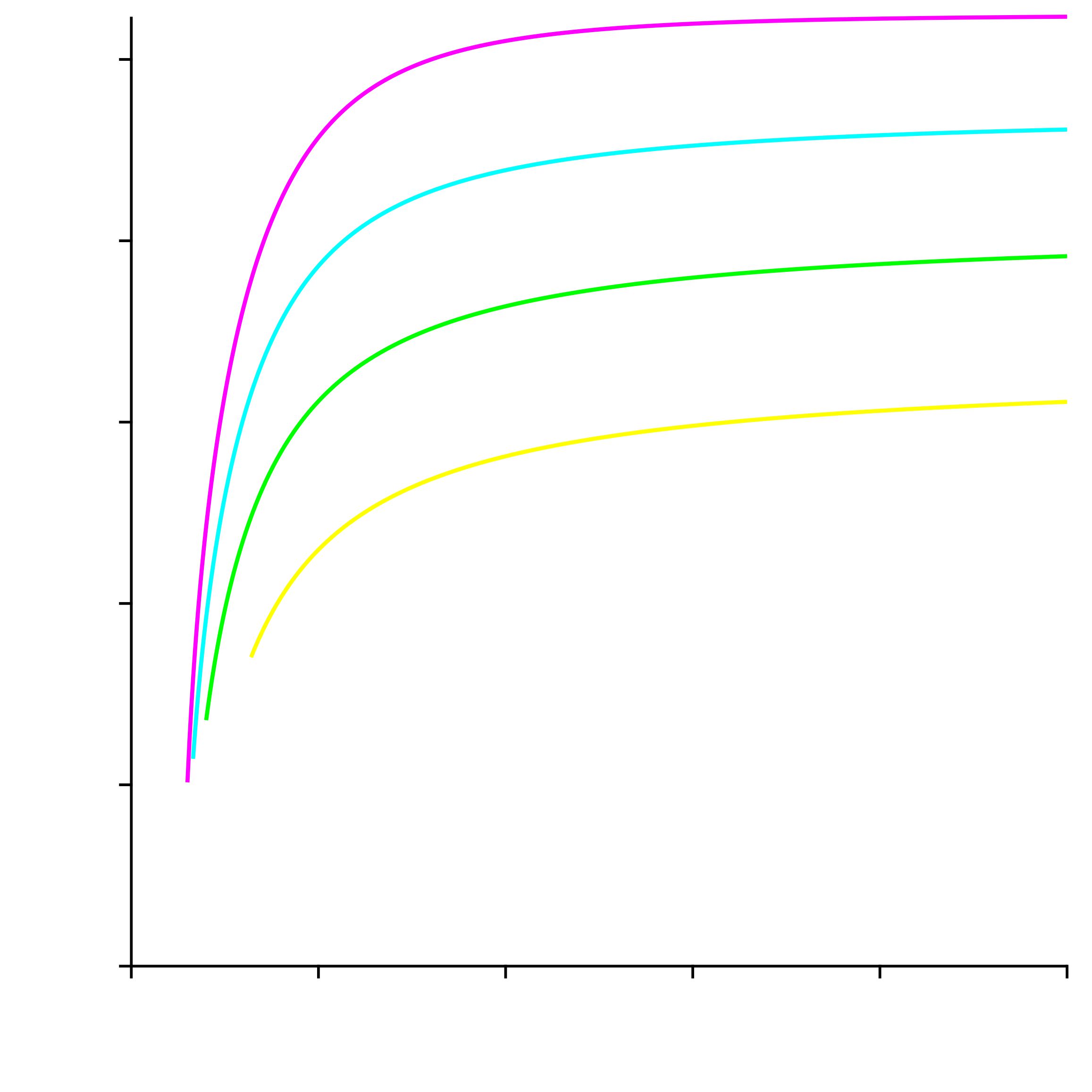}
\begin{small}
\put(-1,90){c)}
\put(9,0){$0$}
\put(26,0){$1$}
\put(-5.5,24.5){$0.2$}
\put(3,8){$0$}
\put(51,1){$Q_1$}
\end{small}     
  \end{overpic}
	\hfill
	$\phm$
\caption{
The change of the height (a), the waist radius (b) and the volume (c) during the snap, where the same color-coding is used as in Fig.\ \ref{fig5}.
}
  \label{fig8}
\end{figure}

\noindent
{\bf Additional graphs.} In the following we give some additional information on 
the families of snapping sandglasses (for $n=3,\ldots ,6$) computed in Section \ref{sec:snap}. 
The change of the height and the waist radius between the closed state and the open one relative to the radius of the $n$-gons in $\alpha$ and $\beta$, respectively, 
is displayed in Fig.\  \ref{fig8}a,b. 
The increase of the volume during the snap relative to the volume of the closed state is plotted in Fig.\  \ref{fig8}c.
This change in volume confirms that there does not exist a continuous isometric deformation between the closed and open state due to 
the {\it Bellows conjecture}.

\end{document}